\documentclass[twoside,11pt,a4paper,leqno]{article}
\usepackage{amssymb,amsmath,amscd,euscript,verbatim,array}
\usepackage[center,pagestyles]{titlesec}
\usepackage{geometry}
\usepackage{anysize}
\usepackage{fancyhdr}
\usepackage{indentfirst}
\usepackage{graphicx}
\usepackage{color}
\usepackage{ifpdf}
\usepackage{mathrsfs}
\usepackage{cite}
\usepackage[numbers,sort&compress]{natbib}
\marginsize{3.5cm}{3.5cm}{2cm}{2cm}

\titleformat{\section}{\centering\normalsize}{\thesection.}{0.5em}{}
\titleformat{\subsection}{\normalsize\bfseries}{\thesubsection.}{0.5em}{}
\titleformat{\subsubsection}{\normalsize\bfseries}{\thesubsubsection.}{0.5em}{}

\newcommand{\R}{\mathbb{R}}

\newtheorem{Theorem}{Theorem}[section]
\newtheorem{Definition}[Theorem]{Definition}
\newtheorem{Lemma}[Theorem]{Lemma}
\newtheorem{Exercise}[Theorem]{Exercise}
\newtheorem{Proposition}[Theorem]{Proposition}
\newtheorem{Remark}[Theorem]{Remark}

\newcommand{\gm}{\gamma}

\newcommand{\A}{\mathbb{A}}

\newcommand{\bthm}{\begin{Theorem}}
\newcommand{\ethm}{\end{Theorem}}
\newcommand{\bpr}{\begin{Proposition}}
\newcommand{\epr}{\end{Proposition}}
\newcommand{\blm}{\begin{Lemma}}
\newcommand{\elm}{\end{Lemma}}
\newcommand{\bex}{\begin{Exercise}}
\newcommand{\eex}{\end{Exercise}}
\newcommand{\be}{\begin{equation}}
\newcommand{\ee}{\end{equation}}
\newcommand{\beal}{\begin{aligned}}
\newcommand{\enal}{\end{aligned}}
\newcommand{\brm}{\begin{Remark}}
\newcommand{\erm}{\end{Remark}}
\newcounter{item}[section]

\newcommand{\Proof}{\textbf{Proof}\hspace{0.3cm}}
\newcommand{\End}{\ensuremath{\hfill{\Box}}\\}
\renewcommand{\title}[1]{\begin{center}\textbf{\large #1}\end{center}}
\renewcommand{\author}[1]{\begin{center}\normalsize #1\end{center}}
\renewcommand{\date}[1]{\begin{center}#1\end{center}}

\makeatletter \@addtoreset{equation}{section}
\makeatother
\begin{document}
\vspace{10pt}
\title{WEAK KAM THEORY FOR GENERAL HAMILTON-JACOBI EQUATIONS II: THE FUNDAMENTAL SOLUTION UNDER  LIPSCHITZ CONDITIONS}

\vspace{6pt}
\author{\sc Lin Wang and Jun Yan}

\vspace{10pt} \thispagestyle{plain}

\begin{quote}
\small {\sc Abstract.} We consider the following evolutionary Hamilton-Jacobi equation with initial condition:
\begin{equation*}
\begin{cases}
\partial_tu(x,t)+H(x,u(x,t),\partial_xu(x,t))=0,\\
u(x,0)=\phi(x),
\end{cases}
\end{equation*}
where $\phi(x)\in C(M,\R)$.
Under some assumptions on the convexity  of  $H(x,u,p)$ with respect to $p$ and the uniform Lipschitz of  $H(x,u,p)$ with respect to $u$, we establish a variational principle and provide  an intrinsic relation between  viscosity solutions and certain minimal characteristics. By introducing an implicitly defined  {\it fundamental solution}, we obtain a variational representation formula of the viscosity solution of the evolutionary  Hamilton-Jacobi equation. Moreover, we discuss the large time behavior of  the viscosity solution of the evolutionary  Hamilton-Jacobi equation  and  provide a dynamical representation formula of the viscosity solution of the stationary Hamilton-Jacobi equation with strictly increasing $H(x,u,p)$ with respect to $u$.
\end{quote}
\begin{quote}
\end{quote}
\begin{quote}
\end{quote} \vspace{25pt}

\tableofcontents
\newpage
\section{\sc Introduction and main results}
Let $M$ be a closed manifold  and $H$ be a $C^r$ ($r\geq 2$) function called a Hamiltonian. We consider the following Hamilton-Jacobi equation:
\begin{equation}\label{hje}
\partial_tu(x,t)+H(x,u(x,t),\partial_xu(x,t))=0,
\end{equation}with the initial condition
\[u(x,0)=\phi(x),\]where $(x,t)\in M\times[0,T]$, $T$ is a positive constant.
The characteristics of (\ref{hje}) satisfies the following equation:
\begin{equation}\label{hjech}
\begin{cases}
\dot{x}=\frac{\partial H}{\partial p},\\
\dot{p}=-\frac{\partial H}{\partial x}-\frac{\partial H}{\partial u}p,\\
\dot{u}=\frac{\partial H}{\partial p}p-H.
\end{cases}
\end{equation}
To avoid the ambiguity, we denote the solution of (\ref{hjech}) (the characteristics of (\ref{hje})) by $(X(t),U(t),P(t))$.

In 1983, M. Crandall and P. L. Lions introduced a notion of weak solution named viscosity solution for overcoming the lack of uniqueness of the solution due to the crossing of characteristics (see \cite{Ar, CL2}). Owing to the notion itself, the uniqueness of the viscosity solution can be followed from comparison principle (see \cite{Ba1,Ba2,CEL,CHL1,CHL2,CL2} for instance).
During the same period, S. Aubry and J. Mather developed a seminar work so called Aubry-Mather theory on global action minimizing orbits for area-preserving twist maps (see \cite{Au,AD,M0,Mat,M2,M3} for instance). Moreover, it was generalized to positive definite Lagrangian systems with multi-degrees of freedom in \cite{M1}.

There is a close connection between viscosity solutions  and Aubry-Mather theory. Roughly speaking, the global minimizing orbits used in Aubry-Mather theory can be embedded into the characteristic fields of PDEs. The similar ideas were reflected  in pioneering papers \cite{E} and \cite{F2} respectively. In \cite{E}, W. E  was concerned with certain weak solutions of Burgers equation. In \cite{F2}, A. Fathi provided a weak solution named weak KAM solution  and implied that the weak KAM solution is a viscosity solution, which initiated so called weak KAM theory.   Later, it was obtained the equivalence between weak KAM solutions and viscosity solutions for the Hamiltonian $H(x,p)$ without the unknown function $u$ under strict convexity and superlinear growth  with respect to $p$.  Moreover, based on the relations between weak KAM solutions and viscosity solutions, the regularity of global subsolutions was improved (see \cite{Be,FS}). A systematic introduction to weak KAM theory can be found in \cite{F3}.

Due to the lack of the variational principle for more general Hamilton-Jacobi equations, the weak KAM theory had been limited to Hamilton-Jacobi equations without the unknown function $u$ explicitly. In \cite{SWY}, the authors made an attempt on the Hamilton-Jacobi equation formed as (\ref{hje}) by a dynamical approach and extended Fathi's weak KAM theory to more general Hamilton-Jacobi equations under the  monotonicity (non-decreasing) and Lipschitz  of $H$ with respect to $u$. Roughly speaking, the weak KAM theory for the Hamilton-Jacobi equations with the unknown function $u$ explicitly is a cornerstone to handle weakly coupled systems and second order equations by a dynamical approach.

In this paper, the  monotonicity (non-decreasing) assumption is dropped, which makes a further step to  enlarge the scope of the weak KAM theory. More precisely, we establish a variational principle and provide  an intrinsic relation between  viscosity solutions and certain minimal characteristics. By introducing an implicitly defined  fundamental solution, we obtain a representation formula of the viscosity solution of (\ref{hje}). Moreover, we discuss the large time behavior of  the viscosity solution of the evolutionary  Hamilton-Jacobi equation  and  provide a dynamical representation formula of the viscosity solution of the stationary Hamilton-Jacobi equation with strictly increasing $H(x,u,p)$ with respect to $u$. Precisely speaking, we are concerned with a  $C^r$ ($r\geq 2$) Hamiltonian $H(x,u,p)$ satisfying the following conditions:
\begin{itemize}
\item [\textbf{(H1)}] \textbf{Positive Definiteness}: $H(x,u,p)$ is strictly convex with respect to $p$;
\item [\textbf{(H2)}] \textbf{Superlinearity in the Fibers}: For every compact set $I$ and any $u\in I$, $H(x,u,p)$ is uniformly superlinear growth with respect  to $p$;
\item [\textbf{(H3)}] \textbf{Completeness of the Flow}: The flows of (\ref{hjech}) generated by $H(x,u,p)$ are complete;
\item [\textbf{(H4)}] \textbf{Uniform Lipschitz}: $H(x,u,p)$ is uniformly Lipschitz with respect to $u$.
\end{itemize}

We use $\mathcal{L}: T^*M\rightarrow TM$ to denote  the Legendre transformation. Let
$\bar{\mathcal{L}}:=(\mathcal{L}, Id)$, where $Id$ denotes the identity map from $\R$ to $\R$. Then $\bar{\mathcal{L}}$ denote a diffeomorphism from $T^*M\times\R$ to $TM\times\R$. By $\bar{\mathcal{L}}$,
the Lagrangian $L(x,u, \dot{x})$ associated to $H(x,u,p)$ can be denoted by
\[L(x,u, \dot{x}):=\sup_p\{\langle \dot{x},p\rangle-H(x,u,p)\}.\]
Let $\Psi_t$ denote the flows of (\ref{hjech}) generated by $H(x,u,p)$. The flows generated by $L(x,u,\dot{x})$
can be denoted by $\Phi_t:=\bar{\mathcal{L}}\circ\Psi_t\circ\bar{\mathcal{L}}^{-1}$. Based on
(H1)-(H4), it follows from $\bar{\mathcal{L}}$ that the Lagrangian $L(x,u,\dot{x})$ satisfies:
\begin{itemize}
\item [\textbf{(L1)}] \textbf{Positive Definiteness}: $L(x,u,\dot{x})$ is strictly convex with respect  to $\dot{x}$;
\item [\textbf{(L2)}] \textbf{Superlinearity in the Fibers}: For every compact set $I$ and any $u\in I$, $L(x,u,\dot{x})$  is uniformly superlinear growth with respect  to $\dot{x}$;
\item [\textbf{(L3)}] \textbf{Completeness of the Flow}: The flows generated by $L(x,u,\dot{x})$  are complete;
\item [\textbf{(L4)}]  \textbf{Uniform Lipschitz}: $L(x,u,\dot{x})$ is uniformly Lipschitz with respect to $u$.
\end{itemize}

If a Hamiltonian $H(x,u,p)$ satisfies (H1)-(H4), then we obtain the following theorem:
\begin{Theorem}\label{two}
For given $x_0, x\in M$, $u_0\in\R$ and $t\in (0, T]$, there exists a unique $h_{x_0,u_0}(x,t)$ satisfying
\begin{equation}
h_{x_0,u_0}(x,t)=u_0+\inf_{\substack{\gm(t)=x \\  \gm(0)=x_0} }\int_0^tL(\gm(\tau),h_{x_0,u_0}(\gm(\tau),\tau),\dot{\gm}(\tau))d\tau,
\end{equation}where the infimums are taken among the absolutely continuous curves $\gm:[0,t]\rightarrow M$. In particular, the infimums are attained at the characteristics of (\ref{hje}).
 Moreover, let $\mathcal{S}_{x_0,u_0}^x$ denote the set of characteristics $(X(t),U(t),P(t))$ satisfying $X(0)=x_0$, $X(t)=x$ and $U(0)=u_0$, then we have
\begin{equation}
h_{x_0,u_0}(x,t)=\inf\left\{U(t):(X(t),U(t),P(t))\in \mathcal{S}_{x_0,u_0}^x\right\}.
\end{equation}
\end{Theorem}

Theorem \ref{two} provides a general variational principle, which builds a bridge between Hamilton-Jacobi equations under (H1)-(H4) and Hamiltonian dynamical systems. As an application, we will obtain a dynamical representation of the viscosity solution of (\ref{hje}).
By analogy with the notion of weak KAM solution of the Hamilton-Jacobi equation without $u$  (see \cite{F3}). We define another weak solution of (\ref{hje}) with initial condition called a variational solution (see Definition \ref{nw}). Based on Theorem \ref{two}, we construct a variational solution of (\ref{hje}) with initial condition.  Following \cite{F3}, we show that the variational solution of (\ref{hje}) is the unique viscosity solution of (\ref{hje}). More precisely, we have the following theorem:
\begin{Theorem}\label{one}
There exists a unique  viscosity solution $u(x,t)$ of (\ref{hje}) with initial condition $u(x,0)=\phi(x)$. Moreover, $u(x,t)$  can be represented as
\begin{equation}
u(x,t)=\inf_{y\in M}h_{y,\phi(y)}(x,t).
\end{equation}
\end{Theorem}

Theorem \ref{two} and Theorem \ref{one} implies the following theorem directly:
\begin{Theorem}\label{three}
For $(x,t)\in M\times [0,T]$, the  viscosity solution $u(x,t)$ of (\ref{hje}) with initial condition $u(x,0)=\phi(x)$ is determined by the minimal characteristic curve. More precisely, we have
\begin{equation}
u(x,t)=\inf_{y\in M}\inf\left\{U(t):(X(t),U(t),P(t))\in \mathcal{S}_{y,\phi(y)}^x\right\},
\end{equation}
where $\mathcal{S}_{y,\phi(y)}^x$ denotes the set of characteristics $(X(t),U(t),P(t))$ satisfying $X(0)=y$, $X(t)=x$ and $U(0)=\phi(y)$.
\end{Theorem}

A similar result corresponding to the viscosity solutions of Hamilton-Jacobi equations without the unknown function $u$ was well known (see Theorem 6.4.6 in \cite{CS} for instance). Theorem \ref{three} implies the relation between the viscosity solutions and  the minimal characteristics still holds for more general Hamilton-Jacobi equations. Roughly speaking, the notion of viscosity solution was invented to avoid the lack of uniqueness owing to the crossing of characteristics.  Based on Theorem \ref{three}, the reason why  the notion of viscosity solution results in the fact without crossing is that the properties of viscosity solutions are determined by certain minimal characteristics.

\begin{Theorem}\label{four}
 There exists an implicitly defined semigroup denoted by $T_t$ such that
\[u(x,t)=T_t\phi(x),\quad h_{y,\phi(y)}(x,s+t)=T_th_{y,\phi(y)}(x,s),\]
where $y\in M,s>0,t\geq 0$ and
\begin{equation}
T_t\phi(x)=\inf_{\gm(t)=x}\left\{\phi(\gm(0))+\int_0^tL(\gm(\tau),T_\tau\phi(\gm(\tau)),
\dot{\gm}(\tau))d\tau\right\}.
\end{equation}
Moreover, for any $\phi(x),\psi(x)\in C(M,\R)$ and $t\in [0,T]$, the solution semigroup $T_t$ has following properties:
\begin{itemize}
\item [I.]  for $\phi\leq\psi$, $T_t\phi\leq T_t\psi$,
\item [II.] $\|T_t\phi-T_t\psi\|_\infty\leq e^{\lambda t}\|\phi-\psi\|_\infty$,
\end{itemize}
where $\lambda>0$ is the Lipschitz constant of $L$.
\end{Theorem}

For $c\in\R$, we denote $L_c:=L+c$. For given $x_0, u_0, x, t$ where $t\in (0,+\infty)$, we define
\begin{equation}
h_{x_0,u_0}^c(x,t)=u_0+\inf_{\substack{\gm(t)=x \\  \gm(0)=x_0} }\int_0^tL_c(\gm(\tau),h_{x_0,u_0}^c(\gm(\tau),\tau),\dot{\gm}(\tau))d\tau,
\end{equation}where the infimums are taken among the absolutely continuous curves $\gm:[0,t]\rightarrow M$.

\begin{Definition}
$c$ is called a critical value if for any $x_0\in M$, $u_0\in\R$ and $t\geq \delta$, it is contained in the following set
\begin{equation}
\mathcal{C}=\left\{c: |h_{x_0,u_0}^c(x,t)|\leq K(u_0)\right\},
\end{equation}
where  $K(u_0)$ is a positive constant depending on $u_0$.
\end{Definition}
 $\mathcal{C}\neq \emptyset$ will be verified in Section 6.
For $a\in\R$, we use $c(L(x,a,\dot{x}))$ to denote Ma\~{n}\'{e} critical value of $L(x,a,\dot{x})$. By \cite{CIPP}, we have
\begin{equation}
c(L(x,a,\dot{x}))=\inf_{u\in C^1(M,\R)}\sup_{x\in M}H(x,a,\partial_xu).
\end{equation}
Without ambiguity, we still use $L$ instead of $L_c$ to denote $L+c$ for $c\in \mathcal{C}$. The same to $H$ and $T_t$. By inspiration of \cite{F22}, the large time behavior of viscosity solutions of Hamilton-Jacobi equations with Hamiltonian independent of $u$ was explored comprehensively  based on both dynamical  and PDE approaches
 (see \cite{ds, II, NR, WY} for instance). Recently, some results on the large time behavior of viscosity solutions of special weakly coupled systems were also obtained (see \cite{CLN,MT1,MT2}). By Theorem \ref{three}, $u(x,t):=T_t\phi(x)$ is the unique viscosity solution of (\ref{hje}) with initial condition $u(x,0)=\phi(x)$.
The following theorem  implies the relation between the viscosity solution of (\ref{hje}) and the one of the stationary equation:
\begin{equation}\label{station}
H(x, u(x), \partial_x u(x))=0.
\end{equation}More precisely, there holds
\begin{Theorem}\label{five}
For any $\phi(x)\in C(M,\R)$, $\liminf_{t\rightarrow\infty}T_t\phi(x)$ exists. Moreover, let \[\underline{u}(x):=\liminf_{t\rightarrow\infty}T_t\phi(x),\]
 then $\underline{u}$ is a weak KAM solution of (\ref{station}).
\end{Theorem}

Let \[h_{x_0,u_0}(x,\infty):=\liminf_{t\rightarrow\infty}h_{x_0,u_0}(x,t),\]
 where $x_0\in M, u_0\in \R$. Based on Theorem \ref{five}, $h_{x_0,u_0}(x,\infty)$ is well-posed. We denote
 \[B(x,u;y):=h_{x,u}(y,\infty)-u.\]

 $B(x,u;y)$ can be referred as the barrier function dented by $h^\infty(x,y)$ in Mather-Fathi theory. Moreover, we define an invariant set called a projected Aubry set as follows;
\begin{equation*}
\mathcal{A}:=\{(x,u)\in M\times\R\ \big|\ B(x,u;x)=0\}.
\end{equation*}
We use $\pi:M\times\R\rightarrow M$ to denote the standard projection via $(x,u)\rightarrow x$.
 \begin{Theorem}\label{six}
Let $H(x,u,p)$ is strictly increasing with respect to $u$ for a given $(x,p)\in T^*M$, then there exists a unique viscosity solution $u(x)$ of (\ref{station}). Moreover,
\begin{equation}
u(x)=\inf_{y\in \pi\mathcal{A}}h_{y,u(y)}(x,\infty).
\end{equation}
\end{Theorem}


\section{\sc Preliminaries}
In this section, we recall the definitions of the weak KAM solution and the viscosity solution of (\ref{hje})  (see \cite{CEL,CL2,F3}). In addition, we provide some aspects of Mather-Fathi theory for  the sake of completeness.
\subsection{Weak KAM solutions and viscosity solutions}
A function $H:TM\rightarrow\R$ called a Tonelli Hamiltonian if $H$ satisfies (H1)-(H2). For the autonomous Hamiltonian systems, the assumption (H3) holds obviously  from the compactness of $M$. The associated Lagrangian is denoted by $L$ via the Legendre transformation. In \cite{F1}, Fathi introduced the definition of the weak KAM solution of negative type of the following Hamilton-Jacobi equation:
\begin{equation}\label{fathi}
H(x,  \partial_x u(x))=0,\quad x\in M,
\end{equation}where $H$ is a Tonelli Hamiltonian.
\begin{Definition}\label{weakkk}
A function $u\in C(M,\R)$ is called a  weak KAM solution of negative type of (\ref{fathi}) if
\begin{itemize}
\item [(i)] for each continuous piecewise $C^1$ curve $\gm:[t_1,t_2]\rightarrow M$ where $t_2>t_1$, we have
\begin{equation}
u(\gm(t_2))-u(\gm(t_1))\leq\int_{t_1}^{t_2}L(\gm(\tau),\dot{\gm}(\tau))d\tau;
\end{equation}
\item [(ii)] for any $x\in M$, there exists a $C^1$ curve $\gm:(-\infty,0]\rightarrow M$ with $\gm(0)=x$ such that for any $t\in (-\infty,0]$, we have
  \begin{equation}
u(x)-u(\gm(t))=\int_{t}^{0}L(\gm(\tau),\dot{\gm}(\tau))d\tau.
\end{equation}
\end{itemize}
\end{Definition}
By  analogy of the definition above, it is easy to define the weak KAM solution of negative type of more general Hamilton-Jacobi equation as follows:
\begin{equation}\label{fathiw}
H(x, u(x), \partial_x u(x))=0,\quad x\in M.
\end{equation}
\begin{Definition}\label{weakkam}
A function $u\in C(M,\R)$ is called a  weak KAM solution of negative type of (\ref{fathiw}) if
\begin{itemize}
\item [(i)] for each continuous piecewise $C^1$ curve $\gm:[t_1,t_2]\rightarrow M$ where $t_2>t_1$, we have
\begin{equation}
u(\gm(t_2))-u(\gm(t_1))\leq\int_{t_1}^{t_2}L(\gm(\tau),u(\gm(\tau)),\dot{\gm}(\tau))d\tau;
\end{equation}
\item [(ii)] for any $x\in M$, there exists a $C^1$ curve $\gm:(-\infty,0]\rightarrow M$ with $\gm(0)=x$ such that for any $t\in (-\infty,0]$, we have
  \begin{equation}
u(x)-u(\gm(t))=\int_{t}^{0}L(\gm(\tau),u(\gm(\tau)),\dot{\gm}(\tau))d\tau.
\end{equation}
\end{itemize}
\end{Definition}

Following from \cite{CEL,CL2,F3},
 a viscosity solution of (\ref{hje}) can be defined as follows:
\begin{Definition}\label{visco}
Let $V$ be an open subset  $V\subset M$,
\begin{itemize}
\item [(i)] A function $u:V\times[0,T]\rightarrow \R$ is a subsolution of  (\ref{hje}), if for every $C^1$ function $\phi:V\times[0,T]\rightarrow\R$ and every point $(x_0,t_0)\in V\times[0,T]$ such that $u-\phi$ has a maximum at $(x_0,t_0)$, we have
\begin{equation}
\partial_t\phi(x_0,t_0)+H(x_0,u(x_0,t_0),\partial_x\phi(x_0,t_0))\leq 0;
\end{equation}
\item [(ii)] A function $u:V\times[0,T]\rightarrow \R$ is a supersolution of  (\ref{hje}), if for every $C^1$ function $\psi:V\times[0,T]\rightarrow\R$ and every point $(x_0,t_0)\in V\times[0,T]$ such that $u-\psi$ has a minimum at $(x_0,t_0)$, we have
\begin{equation}
\partial_t\psi(x_0,t_0)+H(x_0,u(x_0,t_0),\partial_x\psi(x_0,t_0))\geq 0;
\end{equation}
\item [(iii)] A function $u:V\times[0,T]\rightarrow \R$ is a viscosity solution of  (\ref{hje}) on the open subset  $V\subset M$, if it is both a subsolution and a supersolution.
\end{itemize}
\end{Definition}
Under the assumptions (H1)-(H4), it follows from the comparison theorem that the viscosity solution of (\ref{hje}) with initial condition is unique (see \cite{CL2}).

Both of Definition \ref{weakkk} and Definition \ref{weakkam} are concerned with the weak KAM solutions defined on $M\times\R$, while the viscosity solutions of (\ref{hje}) are defined on $M\times[0,T]$. As a bridge connecting them,
we give the definition of another weak solution of (\ref{hje}) with initial condition called a variational solution.
\begin{Definition}\label{nw}
For a given $T>0$, a variational solution of (\ref{hje}) with initial condition is a function $u: M\times [0,T]\rightarrow\R$ for which the following are satisfied:
\begin{itemize}
\item [(i)] for each continuous piecewise $C^1$ curve $\gm:[t_1,t_2]\rightarrow M$ where $0\leq t_1<t_2\leq T$, we have
\begin{equation}
u(\gm(t_2),t_2)-u(\gm(t_1),t_1)\leq\int_{t_1}^{t_2}L(\gm(\tau),u(\gm(\tau),\tau),\dot{\gm}(\tau))d\tau;
\end{equation}
\item [(ii)] for any $0\leq t_1<t_2\leq T$ and $x\in M$, there exists a $C^1$ curve $\gm:[t_1,t_2]\rightarrow M$ with $\gm(t_2)=x$ such that
  \begin{equation}
u(x,t_2)-u(\gm(t_1),t_1)=\int_{t_1}^{t_2}L(\gm(\tau),u(\gm(\tau),\tau),\dot{\gm}(\tau))d\tau.
\end{equation}
\end{itemize}
\end{Definition}
The existence of  the variational solutions will be verified in Section 4.

\subsection{The minimal action and the  fundamental solution}
Let $L:TM\rightarrow\R$ be a Tonelli Lagrangian. We define the function $h_t:M\times M\rightarrow\R$ by
\begin{equation}\label{mather}
h_t(x,y)=\inf_{\substack{\gm(0)=x\\ \gm(t)=y}}\int_0^t L(\gm(\tau),\dot{\gm}(\tau))d\tau,
\end{equation}
where the infimums are taken among the absolutely continuous curves $\gm:[0,t]\rightarrow M$.
By Tonelli theorem (see \cite{F3,M1}), the infimums in (\ref{mather}) can be achived. Let $\bar{\gm}$ be an absolutely continuous curve with $\bar{\gm}(0)=x$ and $\bar{\gm}(t)=y$ such that the infinmum is achieved at $\bar{\gm}$. Then $\bar{\gm}$ is called a minimal curve. By \cite{M1}, the minimal curves satisfy the Euler-Lagrange equation generated by $L$. The quantity $h_t(x,y)$ is called a minimal  action.
From the definition of $h_t(x,y)$, it follows that for each $x,y,z\in M$ and each $t,t'>0$, we have
 \begin{equation}\label{yan}
h_{t+t'}(x,z)\leq h_{t}(x,y)+h_{t'}(y,z).
\end{equation}
  In particular,  we have
\begin{equation}\label{wang}
h_{t+t'}(x,y)=h_{t}(x,\bar{\gm}(t))+h_{t'}(\bar{\gm}(t),y),
\end{equation}
where $\bar{\gm}$ is  a minimal curve with $\bar{\gm}(0)=x$ and $\bar{\gm}(t+t')=y$.

Consider the following Hamilton-Jacobi equation:
 \begin{equation}\label{fathitt}
\begin{cases}
\partial_tu(x,t)+H(x,  \partial_xu(x,t))=0,\\
u(x,0)=\phi(x),
\end{cases}
\end{equation}
where $\phi(x)\in C(M)$.
By \cite{F3}, a viscosity solution of (\ref{fathitt}) can be represented as
 \begin{equation}\label{infc}
u(x,t):=\inf_{y\in M}\left\{\phi(y)+h^t(y,x)\right\}.
\end{equation}
The right side of (\ref{infc}) is also called inf-convolution of $\phi$, due to the formal analogy with the usual convolution (see \cite{CS}). Moreover, the minimal action $h^t(y,x)$ can be viewed as a {\it fundamental solution} of (\ref{fathitt}) (see \cite{GT}).

The following conception is crucial in our context.
\begin{Definition}\label{cali}
For $u(x,t)\in C(M\times[0,T],\R)$, a curve $\gm:I\rightarrow M$ is called a calibrated curve of $u$ if for every $t_1, t_2\in I$ with $0\leq t_1< t_2$, we have
\[u(\gm(t_2),t_2)=u(\gm(t_1),t_1)+\int_{t_1}^{t_2}L(\gm(\tau),u(\gm(\tau),\tau),\dot{\gm}(\tau))d\tau.\]
\end{Definition}

We are devoted to detecting the viscosity solution of (\ref{hje}) from a dynamical view.  For given $x_0, x\in M$, $u_0\in\R$ and $t\in (0, T]$, we define formally:
\begin{equation}\label{u}
h_{x_0,u_0}(x,t)=u_0+\inf_{\substack{\gm(t)=x \\  \gm(0)=x_0} }\int_0^tL(\gm(\tau),h_{x_0,u_0}(\gm(\tau),\tau),\dot{\gm}(\tau))d\tau,
\end{equation}where the infimums are taken among the absolutely continuous curves $\gm:[0,t]\rightarrow M$. It is easy to see that the cure achieving the infimum in the right side of (\ref{u}) is a calibrated curve of $h_{x_0,u_0}(x,t)$. To fix the notions, we call $h_{x_0,u_0}(x,t)$ the fundamental solution of (\ref{hje}).
  In next section, we will show the well-posedness of $h_{x_0,u_0}(x,t)$ under the assumptions (L1)-(L4).

\section{\sc Variational principle}
In this section, we are devoted to proving Theorem \ref{two}.
The proof will be proceeded  by four steps. In the first step, we will prove the existence and uniqueness  of $h_{x_0,u_0}(x,t)$. In the second step, we will verify $h_{x_0,u_0}(x,t)$ to satisfy a triangle inequality. In the third step, we will show that the relation between calibrated curves and characteristics. Based on the preliminaries in former  steps, we will give a relation between $h_{x_0,u_0}(x,t)$ and $U(t)$ belonging to  a characteristic curve $(X(t),U(t),P(t))$ in the last step.

\subsection{Existence and uniqueness of the fundamental solution}
In this step, we are concerned with the existence and uniqueness  of $h_{x_0,u_0}(x,t)$.
We use $C^{\text{ac}}([0,t],M)$ to denote the set of all absolutely continuous curves $\gm:[0,t]\rightarrow M$. First of all, we verify the existence of $h_{x_0,u_0}(x,t)$.
\begin{Lemma}\label{exist}
There exists $h_{x_0,u_0}(x,t)\in C(M\times(0,T],\R)$ such that
\begin{equation}
h_{x_0,u_0}(x,t)=u_0+\inf_{\substack{\gm(t)=x \\  \gm(0)=x_0 } }\int_0^tL(\gm(\tau),h_{x_0,u_0}(\gm(\tau),\tau),\dot{\gm}(\tau))d\tau,
\end{equation}where $\gm\in C^{\text{ac}}([0,t],M)$.
\end{Lemma}
\Proof  For the simplicity of notations, without ambiguity, we drop the subscripts $x_0$ and $u_0$ of $h_{x_0,u_0}(x,t)$. We consider a sequence generated by the following iteration:
\begin{equation}\label{ui}
h_{i+1}(x,t)=u_0+\inf_{ \gm(t)=x \atop  \gm(0)=x_0 }\int_0^tL(\gm(\tau),h_{i}(\gm(\tau),\tau),\dot{\gm}(\tau))d\tau,
\end{equation}where $i=0,1,2,\ldots$ and $h_0(x,t)=u_0$.
By means of a simple modification of Tonelli's theorem (see \cite{F3} and \cite{M1}), we have that for a given $h_i(x,t)\in C(M\times(0,T],\R)$ there exists an absolutely continuous curve $\gm_i: [0,t]\rightarrow M$ satisfying $\gm_i(0)=x_0$  and $\gm_i(t)=x$ such that the infimum in (\ref{u}) can be achieved. To fix the notions, $\gm_i$ is called a minimal curve of $h_i$.

Let $\bar{\gm}:[0,t]\rightarrow M$ be an absolutely continuous curve satisfying $\bar{\gm}(0)=x_0$ and $\bar{\gm}(t)=x$. By the construction of $h_i$, there holds for $s\in [0,t]$,
\begin{equation}
h_{1}(\gm(s),s)=u_0+\inf_{ \gm(s)=\bar{\gm}(s) \atop  \gm(0)=x_0 }\int_0^sL(\gm(\tau),u_0,\dot{\gm}(\tau))d\tau,
\end{equation}
It is easy to see that
\begin{equation}\label{bound}
|h_{1}(\gm(s),s)-u_0|\leq K,
\end{equation}
where $K$ is a positive constant independent of $s$.  Let $\gm_2:[0,t]\rightarrow M$ be a minimal curve of $h_2$ with $\gm_2(0)=x_0$ and  $\gm_2(t)=x$. Let $\gm_1:[0,t]\rightarrow M$ be a minimal curve of $h_1$ with $\gm_1(0)=x_0$ and  $\gm_1(s)=\gm_2(s)$. By (L4), we have
\begin{equation}\label{lip}
|L(x,u,\dot{x})-L(x,v,\dot{x})|\leq \lambda|u-v|.
\end{equation}Then for $s\in [0,t]$, we have
\begin{align*}
&\ \ \ h_2(\gm_2(s),s)-h_1(\gm_2(s),s)\\
&\leq\int_0^s L(\gm_1(\tau),h_1(\gm_1(\tau),\tau),\dot{\gm}_1(\tau))d\tau-\int_0^sL(\gm_1(\tau),u_0,\dot{\gm}_1(\tau))d\tau,\\
&\leq\int_0^s |L(\gm_1(\tau),h_1(\gm_1(\tau),\tau),\dot{\gm}_1(\tau))-L(\gm_1(\tau),u_0,\dot{\gm}_1(\tau))|d\tau,\\
&\leq\lambda\int_0^s|h_1(\gm_1(\tau),\tau)-u_0|d\tau,
\end{align*}
which together with (\ref{bound}) implies
\begin{equation}\label{2kt}
h_2(\gm_2(s),s)-h_1(\gm_2(s),s)\leq C\lambda s,
\end{equation} By a similar argument, we have
\begin{equation}
 h_2(\gm_2(s),s)-h_1(\gm_2(s),s)\geq -C\lambda s.
\end{equation}
In particular, there holds for  $(x,t)\in M\times (0,T]$,
\[|h_2(x,t)-h_1(x,t)|\leq C\lambda t.\]
 Let $\gm_3:[0,t]\rightarrow M$ be a minimal curve of $h_3$ with $\gm_3(0)=x_0$ and  $\gm_3(t)=x$. Let $\gm_2:[0,t]\rightarrow M$ be a minimal curve of $h_2$ with $\gm_2(0)=x_0$ and  $\gm_2(s)=\gm_3(s)$. Moreover,
\begin{align*}
&\ \ \ h_3(\gm_3(s),s)-h_2(\gm_3(s),s)\\
&\leq\int_0^s L(\gm_2(\tau),h_2(\gm_2(\tau),\tau),\dot{\gm}_2(\tau))d\tau-\int_0^s L(\gm_2(\tau),h_1(\gm_2(\tau),\tau),\dot{\gm}_2(\tau))d\tau,\\
&\leq\int_0^s |L(\gm_2(\tau),h_2(\gm_2(\tau),\tau),\dot{\gm}_2(\tau))-L(\gm_2(\tau),h_1(\gm_2(\tau),\tau),\dot{\gm}_2(\tau))|d\tau,\\
&\leq\lambda\int_0^s |h_2(\gm_2(\tau),\tau))-h_1(\gm_2(\tau),\tau)|d\tau\leq \lambda^2 C\int_0^s \tau d\tau=\frac{1}{2}C(\lambda s)^2.
\end{align*}
 By a similar argument, we have
 \begin{equation}
h_3(\gm_3(s),s)-h_2(\gm_3(s),s)\geq-\frac{1}{2}C(\lambda s)^2.
\end{equation}
In particular, we have
\begin{equation}\label{C32}
|h_3(x,t)-h_2(x,t)|\leq\frac{1}{2}C(\lambda t)^2,
\end{equation}
Repeating the argument above $n$ times, we have
\begin{equation}\label{Cn}
|h_{n+1}(x,t)-h_n(x,t)|\leq \frac{1}{n!}C(\lambda t)^n.
\end{equation}It follows from (\ref{Cn}) that as $n\rightarrow \infty$,
\begin{equation}
|h_{n+1}(x,t)-h_n(x,t)|\rightarrow 0,
\end{equation}which implies that $\{h_n\}$ is a Cauchy sequence, hence there exists $\bar{h}(x,t)\in C(M\times (0,T],\R)$ such that
\begin{equation}
\lim_{n\rightarrow\infty}h_n(x,t)=\bar{h}(x,t),
\end{equation}where $\bar{h}(x,t)$ satisfies (\ref{u}). This finishes the proof of Lemma \ref{exist}.
\End

Lemma \ref{exist} implies that there exists $h_{x_0,u_0}(x,t)\in C(M\times(0,T],\R)$ such that
\begin{equation}\label{hh}
h_{x_0,u_0}(x,t)=u_0+\inf_{\substack{\gm(t)=x \\  \gm(0)=x_0 } }\int_0^tL(\gm(\tau),h_{x_0,u_0}(\gm(\tau),\tau),\dot{\gm}(\tau))d\tau.
\end{equation} In particular, the infimum can be achieved at an absolutely continuous curve denoted by $\bar{\gm}$. By Definition \ref{cali}, $\bar{\gm}$ is a calibrated curve.

The following lemma implies the uniqueness of $h_{x_0,u_0}(x,t)$.

\begin{Lemma}\label{unique}
If both $h_{x_0,u_0}(x,t)$ and $g_{x_0,u_0}(x,t)$ satisfy (\ref{u}), then $h_{x_0,u_0}(x,t)=g_{x_0,u_0}(x,t)$ for $(x,t)\in M\times (0,T]$.
\end{Lemma}
The proof of Lemma \ref{unique} depends on a useful inequality as follows.

\textbf{Gronwall's inequality:} Let $F:[0,t]\rightarrow\R$ be continuous and nonnegative. Suppose $C\geq 0$ and $\lambda\geq 0$ are such that for any $s\in [0,t]$,
 \begin{equation}
F(s)\leq C+\int_{0}^s\lambda F(\tau)d\tau.
\end{equation}Then, for any $s\in [0,t]$,
 \begin{equation}
F(s)\leq Ce^{\lambda s}.
\end{equation}
Taking $C=0$, we have the following lemma.
\begin{Lemma}\label{gron}
Let $F:[0,t]\rightarrow\R$ be continuous with $F(0)=0$ and $F(s)>0$ for $s\in (0,t]$. Then for a given $\lambda\geq 0$, there exists $s_0\in (0,t]$ such that
 \begin{equation}\label{Ft}
F(s_0)>\int_{0}^{s_0}\lambda F(\tau)d\tau.
\end{equation}
\end{Lemma}

\textbf{Proof of Lemma \ref{unique}:} The same as the notations in Lemma \ref{exist}, we denote $h_{x_0,u_0}(x,t)$ and $g_{x_0,u_0}(x,t)$ by $h(x,t)$ and $g(x,t)$ respectively.

On one hand, we will prove
\begin{equation}
h(x,t)\leq g(x,t).
\end{equation}
By contradiction, we assume $h(x,t)> g(x,t)$.
 Let $\gm_g$ be a calibrated curve of $g$ with $\gm_g(0)=x_0$, $\gm_g(t)=x$. We denote
 \begin{equation}
F(\tau)=h(\gm_g(\tau),\tau)-g(\gm_g(\tau),\tau),
\end{equation}where $\tau\in [0,t]$. By (\ref{u}), we have $F(0)=0$. The assumption $h(x,t)> g(x,t)$ implies $F(t)>0$. Hence, there exists $\tau_0\in [0,t)$ such that $F(\tau_0)=0$ and  $F(\tau)> 0$ for $\tau> \tau_0$.
Let $\gm_h$ be a calibrated curve of $h$ with $\gm_h(0)=x_0$, $\gm_h(\tau_0)=\gm_g(\tau_0)$. For $s\in [\tau_0,t]$, we construct $\gm_s:[0,s]\rightarrow M$ as follows:
\begin{equation}
\gm_s(\tau)=\left\{\begin{array}{ll}
\hspace{-0.4em}\gm_h(\tau),& \tau\in [0,\tau_0],\\
\hspace{-0.4em}\gm_g(\tau),&\tau\in (\tau_0,s].\\
\end{array}\right.
\end{equation}
 Based on the definition of $h(x,t)$ (see (\ref{u})), we have
 \begin{align*}
 h(\gm_g(s),s)&=u_0+\inf_{\substack{\gm(s)=\gm_g(s) \\  \gm(0)=x_0 } }\int_0^sL(\gm(\tau),h(\gm(\tau),\tau),\dot{\gm}(\tau))d\tau,\\
 &\leq u_0+\int_{0}^sL(\gm_s(\tau),h(\gm_s(\tau),\tau),\dot{\gm}_s(\tau))d\tau,\\
 &= h(\gm_g(\tau_0),\tau_0)+\int_{\tau_0}^sL(\gm_g(\tau),h(\gm_g(\tau),\tau),\dot{\gm}_g(\tau))d\tau.
 \end{align*}
Similarly, for $g(x,t)$, we have
 \begin{equation}
g(\gm_g(s),s)=g(\gm_g(\tau_0),\tau_0)+\int_{\tau_0}^sL(\gm_g(\tau),g(\gm_g(\tau),\tau),\dot{\gm}_g(\tau))d\tau.
\end{equation}
Since $h(\gm_g(\tau_0),\tau_0)-g(\gm_g(\tau_0),\tau_0)=F(\tau_0)=0$,   then we have
 \begin{align*}
 &h(\gm_g(s),s)-g(\gm_g(s),s)\\
 &\leq \int_{\tau_0}^sL(\gm_g(\tau),h(\gm_g(\tau),\tau),\dot{\gm}_g(\tau))-
 L(\gm_g(\tau),g(\gm_g(\tau),\tau),\dot{\gm}_g(\tau))d\tau,\\
 &\leq \int_{\tau_0}^s\lambda |h(\gm_g(\tau),\tau)-g(\gm_g(\tau),\tau)|d\tau,\\
 &=\int_{\tau_0}^s\lambda (h(\gm_g(\tau),\tau)-g(\gm_g(\tau),\tau))d\tau,
 \end{align*}where the second inequality is owing to (A2).
It follows that for any $s\in (\tau_0,t]$
 \begin{equation}\label{Ft1}
F(s)\leq \int_{\tau_0}^s\lambda F(\tau)d\tau,
\end{equation}
which is in contradiction with Lemma \ref{gron}.  Thus, we obtain $h(x,t)\leq g(x,t)$.

On the other hand, it follows from a similar argument that $h(x,t)\geq g(x,t)$. So far, we have shown that
$h(x,t)=g(x,t)$ for $(x,t)\in M\times (0,T]$, which finishes the proof of Lemma \ref{unique}.
\End

\subsection{A triangle inequality}
 Lemma \ref{exist} and Lemma \ref{unique} imply the well definiteness of $h_{x_0,u_0}(x,t)$. For the simplicity of notations,  we drop the subscripts $x_0$ and $u_0$ of $h_{x_0,u_0}(x,t)$.  The following lemma implies that $h(x,t)$ satisfies a triangle inequality.

 \begin{Lemma}\label{tria1}
  \begin{equation}
h(x,t+s)=\inf_{y\in M}h_{y,h(y,t)}(x,s).
\end{equation}
 \end{Lemma}
\Proof
On one hand, we will prove $h(x,t+s)\geq\inf_{y\in M}h_{y,h(y,t)}(x,s)$. Let $\gm_1:[0,t+s]\rightarrow M$ be a calibrated curve of $h$ with $\gm_1(0)=x_0$ and $\gm_1(t+s)=x$. Consider $\bar{y}\in \gm_1$ with $\gm_1(t)=\bar{y}$. It suffices to show
  \begin{equation}\label{s11}
h(x,t+s)\geq h_{\bar{y},h(\bar{y},t)}(x,s).
\end{equation}
By contradiction, we assume $h(x,t+s)< h_{\bar{y},h(\bar{y},t)}(x,s)$. By the definition of $h(x,t+s)$, we have
\begin{align*}
h(x,t+s)&=u_0+\int_0^{t+s}L(\gm_1(\tau),h(\gm_1(\tau),\tau),\dot{\gm}_1(\tau))d\tau,\\
&=h(\bar{y},t)+\int_t^{t+s}L(\gm_1(\tau),h(\gm_1(\tau),\tau),\dot{\gm}_1(\tau))d\tau,\\
&=h(\bar{y},t)+\int_0^{s}L(\gm_1(\sigma+t),h(\gm_1(\sigma+t),\sigma+t),\dot{\gm}_1(\sigma+t))d\sigma.
\end{align*}
By the definition of $h_{\bar{y},h(\bar{y},t)}(x,s)$, we have
\begin{align*}
h_{\bar{y},h(\bar{y},t)}(x,s)&=h(\bar{y},t)+\inf_{\substack{\gm(s)=x \\  \gm(0)=\bar{y} } }
\int_0^{s}L(\gm(\tau),h_{\bar{y},h(\bar{y},t)}(\gm(\tau),\tau),\dot{\gm}(\tau))d\tau,\\
&\leq h(\bar{y},t)+\int_0^{s}L(\gm_1(\sigma+t),h_{\bar{y},h(\bar{y},t)}(\gm_1(\sigma+t),\sigma),\dot{\gm}(\sigma+t))d\sigma.
\end{align*}
For the simplicity of notations, we denote $u(\sigma):=h(\gm_1(\sigma+t),\sigma+t)$ and $v(\sigma):=h_{\bar{y},h(\bar{y},t)}(\gm_1(\sigma+t),\sigma)$. In particular, we have $h(x,t+s)=u(s)$ and $h_{\bar{y},h(\bar{y},t)}(x,s)=v(s)$. Let
 \begin{equation}
F(\sigma):=v(\sigma)-u(\sigma),
\end{equation}where $\sigma\in [0,s]$. It is easy to see that $u(0)=h(\bar{y},t)=v(0)$. Then we have $F(0)=0$. The assumption $h(x,t+s)< h_{\bar{y},h(\bar{y},t)}(x,s)$ implies $F(s)>0$. Hence, there exists $\sigma_0\in [0,s)$ such that $F(\sigma_0)=0$ and  $F(\sigma)> 0$ for $\sigma> \sigma_0$. Moreover, for any $\tau\in (\sigma_0,s]$, we have
  \begin{equation}
u(\tau)=u(\sigma_0)+\int_{\sigma_0}^{\tau}L(\gm_1(\sigma+t),u(\sigma),\dot{\gm}_1(\sigma+t))d\sigma.
\end{equation}
Let $\gm_2$ be a calibrated curve of $h$ with $\gm_2(0)=\bar{y}$, $\gm_2(\sigma_0)=\gm_1(\sigma_0+t)$. For $\sigma\in [\sigma_0,\tau]$, we construct $\gm_\tau:[0,\tau]\rightarrow M$ as follows:
\begin{equation}
\gm_\tau(\sigma)=\left\{\begin{array}{ll}
\hspace{-0.4em}\gm_2(\sigma),& \sigma\in [0,\sigma_0],\\
\hspace{-0.4em}\gm_1(\sigma+t),&\sigma\in (\sigma_0,\tau].\\
\end{array}\right.
\end{equation} Moreover,
for any $\tau\in (0,s]$, we have
  \begin{align*}
v(\tau)&\leq h(\bar{y},t)+\int_0^{\tau}L(\gm_\tau(\sigma+t),
h_{\bar{y},h(\bar{y},t)}(\gm_\tau(\sigma),\sigma),\dot{\gm}_\tau(\sigma))d\sigma,\\
&=v(\sigma_0)+\int_{\sigma_0}^{\tau}L(\gm_1(\sigma+t),v(\sigma),\dot{\gm}_1(\sigma+t))d\sigma.
\end{align*}
Since $v(\sigma_0)-u(\sigma_0)=F(\sigma_0)=0$, it follows from (A2) that
  \begin{equation}
v(\tau)-u(\tau)\leq \int_{\sigma_0}^\tau\lambda(v(\sigma)-u(\sigma))d\sigma.
\end{equation}
It yields that for any $\tau\in (\sigma_0,s]$,
 \begin{equation}
F(\tau)\leq \int_{\sigma_0}^\tau\lambda F(\sigma)d\sigma,
\end{equation}
which is in contradiction with Lemma \ref{gron}. Hence,
 \begin{equation}\label{s22}
h(x,t+s)\geq\inf_{y\in M}h_{y,h(y,t)}(x,s).
\end{equation}

On the other hand, we will prove $h(x,t+s)\leq\inf_{y\in M}h_{y,h(y,t)}(x,s)$. Due to the compactness of $M$, there exists $\tilde{y}$ such that $\inf_{y\in M}h_{y,h(y,t)}(x,s)=h_{\tilde{y},h(\tilde{y},t)}(x,s)$.
 Let $\gm_3:[0,t]\rightarrow M$ be a calibrated curve of $h$ with $\gm_3(0)=x_0$ and $\gm_3(t)=\tilde{y}$.  Let $\gm_4:[0,s]\rightarrow M$ be a calibrated curve of $h_{\tilde{y},h(\tilde{y},t)}$ with $\gm_4(0)=\tilde{y}$ and $\gm_4(s)=x$. By a similar argument as (\ref{s11})-(\ref{s22}), we have
 \begin{equation}\label{s23}
h(x,t+s)\leq\inf_{y\in M}h_{y,h(y,t)}(x,s),
\end{equation}
which together with (\ref{s22}) implies
  \begin{equation}
h(x,t+s)=\inf_{y\in M}h_{y,h(y,t)}(x,s).
\end{equation}
This completes the proof of Lemma \ref{tria1}.
\End
From Lemma \ref{tria1}, we obtain
  \begin{equation}\label{triin}
h(x,t+s)-u_0\leq (h_{y,h(y,t)}(x,s)-h(y,t))+(h(y,t)-u_0),
\end{equation}
which can be seen as a triangle inequality. In particular, the equality holds if and only if $y$ belongs to the calibrated curve $\gm$ of $h$ with $\gm(0)=x_0$ and $\gm(t+s)=x$.
\subsection{Calibrated curves and characteristics}
In Step 1, we obtain that there exists a unique  $h(x,t)\in C(M\times(0,T],\R)$ satisfying (\ref{hh}) and a calibrated curve $\gm$ of $h$. In this step, we will show that the relation between calibrated curves and characteristics. More precisely, we have the following lemma:

\begin{Lemma}\label{chc}
Let $\bar{\gm}:[0,t]\rightarrow M$ be a calibrated curve of $h$, then $\bar{\gm}$ is $C^1$ and for $\tau\in (0,t)$,
$(\bar{\gm}(\tau),u(\tau),p(\tau))$ satisfies the characteristics equation (\ref{hjech}) where
\begin{equation}\label{ccure}
u(\tau)=h(\bar{\gm}(\tau),\tau)\quad\text{and}\quad p(\tau)=\frac{\partial L}{\partial \dot{x}}(\bar{\gm}(\tau),h(\bar{\gm}(\tau),\tau),\dot{\bar{\gm}}(\tau)).
\end{equation}
\end{Lemma}
\Proof
Since $\bar{\gm}\in C^{ac}([0,t],M)$, then the derivative $\dot{\bar{\gm}}(\tau)$ exists almost everywhere for $\tau\in [0,t]$. Let $t_0\in (0,t)$ be a differentiate point of $\bar{\gm}(\tau)$. Without loss of generality, we assume $0<t_0<t$. For the simplicity of notations and without ambiguity, we denote
\begin{equation}
(x_0,u_0,v_0):=(\bar{\gm}(t_0),h(\bar{\gm}(t_0),t_0),\dot{\bar{\gm}}(t_0)).
 \end{equation}

First of all, we will construct a classical solution on a cone-like region (see (\ref{tri}) below). Let $k:=|v_0|$ and
 \[B(0,2k):=\{v:|v|< 2k,\ v\in T_{x_0}M\}.\]
We use $B^*(0,2k)$ to denote the image of $B(0,2k)$ via the Legendre transformation $\mathcal{L}^{-1}:TM\rightarrow T^*M$. That is
\[B^*(0,2k):=\left\{p: p=\frac{\partial L}{\partial v}(x_0,u_0,v),\ v\in B(0,2k)\right\}.\]
  Let $\Psi_t:T^*M\times\R\rightarrow T^*M\times\R$ denote the follow generated by the characteristics equation (\ref{hjech}). Let $\pi$ be a projection from $T^*M\times\R$ to $ T^*M$ via $(x,p,u)\rightarrow (x,p)$ and let
$B_t^*(0,2k):=\pi\circ\Psi_{t-t_0}(B^*(0,2k),u_0)$. We denote
\[\Pi_t:B_t^*(0,2k)\rightarrow M.\]
Since the Legendre transformation $\mathcal{L}$ is a diffeomorphism, then for given   $\epsilon>0$ small enough and $\tau\in [t_0,t_0+\epsilon]$,  $\Pi_\tau$ is a diffeomorphism onto the image denoted by $\Omega_\tau:=\Pi_\tau( B_\tau^*(0,2k))$.
 We use $\Omega^\epsilon$ to denote the following cone-like region:
\begin{equation}\label{tri}
\Omega^\epsilon:=\{(\tau,x): \tau\in (t_0,t_0+\epsilon),\ x\in \Omega_\tau\}.
 \end{equation}Then for any $(\tau,x)\in\Omega^\epsilon$, there exists a unique $p_0\in B^*(0,2k)$ such that $X(\tau)=x$ where
 \[(X(\tau), U(\tau),P(\tau)):=\Phi_\tau(x_0,u_0,p_0).\]
 Hence, for  any $(\tau,x)\in \Omega^\epsilon$, one can define a $C^1$ function by $S(x,\tau)=U(\tau)$. In particular, we have $S(x,t_0)=u_0$. Moreover, it follows from  the method of characteristics (see \cite{F3,Li} for instance) that $S(x,\tau)$
is a solution of the following equation:
\begin{equation}\label{sss}
\partial_\tau S(x,\tau)+H(x,S(x,\tau),\partial_x S(x,\tau))=0
\end{equation}with $\partial_x S(x_0,t_0)=\partial_v L(x_0,S(x_0,t_0),v_0)$, where $L$ denotes Lagrangian via the Legendre transformation associated to the Hamiltonian $H$. Fix $\tau\in [t_0,t_0+\epsilon]$ and let $S_\tau(x):=S(x,\tau)$. We denote
\begin{equation}
\text{grad}_LS_\tau(x):=\frac{\partial H}{\partial p}(x,S_\tau(x),p),
\end{equation}where $p=\partial_x S_\tau(x)$. In particular, we have $v_0=\text{grad}_LS_{t_0}(x_0)$. It is easy to see that $\text{grad}_LS_\tau(x)$ gives rise to a vector field on $M$. Let $\Omega$ be the Legendre transformation of $\Omega^*$. Moreover, we have the following claim:

\textbf{Claim:} Let $\gm$ be an absolutely continuous curve with $(\tau,\gm(\tau))\in \Omega$ for  $\tau\in [a,b]\subset [t_0,t_0+\epsilon]$, we have
\begin{equation}\label{hjee}
S(\gm(b),b)-S(\gm(a),a)\leq\int_a^bL(\gm(\tau),
S(\gm(\tau),\tau),\dot{\gm}(\tau))d\tau,
\end{equation}
where the equality holds if and only if $\gm$ is a trajectory of the vector field $\text{grad}_LS_\tau(x)$.

\Proof
From the regularity of $S(x,\tau)$, it follows that
\begin{equation}\label{sasb}
S(\gm(b),b)-S(\gm(a),a)=\int_a^b\left\{\frac{\partial S}{\partial t}(\gm(\tau),\tau)+\langle \frac{\partial S}{\partial x}(\gm(\tau),\tau),\dot{\gm}(\tau)\rangle\right\}d\tau.
\end{equation}By virtue of Fenchel inequality, for each $\tau$ where $\dot{\gm}(\tau)$ exists, we have
\begin{align*}
\langle \frac{\partial S}{\partial x}(\gm(\tau),\tau),\dot{\gm}(\tau)\rangle\leq &H(\gm(\tau),S(\gm(\tau),\tau),\frac{\partial S}{\partial x}(\gm(\tau),\tau))\\
&+L(\gm(\tau),
S(\gm(\tau),\tau),\dot{\gm}(\tau)).
\end{align*}It follows from (\ref{sss}) that for almost every $\tau\in [a,b]$
\begin{equation}
\frac{\partial S}{\partial t}(\gm(\tau),\tau)+\langle \frac{\partial S}{\partial x}(\gm(\tau),\tau),\dot{\gm}(\tau)\rangle\leq L(\gm(\tau),
S(\gm(\tau),\tau),\dot{\gm}(\tau)).
\end{equation}
By integration, it follows from (\ref{sasb})that
\begin{equation}\label{ss}
S(\gm(b),b)-S(\gm(a),a)\leq\int_a^bL(\gm(\tau),
S(\gm(\tau),\tau),\dot{\gm}(\tau))d\tau.
\end{equation}We have equality in (\ref{ss}) if and only if the equality holds in the Fenchel inequality, i.e. $\dot{\gm}(\tau)=\text{grad}_LS_\tau(x)$ which means that $\gm$ is a trajectory of the vector field $\text{grad}_LS_\tau(x)$.
\End

Based on the construction of $\Omega^\epsilon$, we have $(\tau,\bar{\gm}(\tau))\in \Omega^\epsilon$ for $\tau\in [t_0,t_0+\epsilon]$. By virtue of a similar argument as the one in the proof of Lemma \ref{unique}, we have for any $\tau\in [t_0,t_0+\epsilon]$,
\begin{equation}
S(\bar{\gm}(\tau),\tau)=h(\bar{\gm}(\tau),\tau),
\end{equation}where $\bar{\gm}$ is a calibrated curve of $h$ with $\bar{\gm}(t_0)=x_0$. From the definition of $h_{x_0,u_0}$ (see (\ref{u})), it follows that
\begin{align*}
S(\bar{\gm}(t_0+\epsilon),t_0+\epsilon)=S(\bar{\gm}(t_0),t_0)
+\int_{t_0}^{t_0+\epsilon}L(\bar{\gm}(\tau),S(\bar{\gm}(\tau),\tau),\dot{\bar{\gm}}(\tau))d\tau,
\end{align*}which implies $\bar{\gm}(\tau)$  is a trajectory of the vector field $\text{grad}_LS_\tau(x)$. Let
\begin{equation}
u(\tau):=h(\bar{\gm}(\tau),\tau)\quad\text{and}\quad p(\tau):=\frac{\partial L}{\partial \dot{x}}(\bar{\gm}(\tau),h(\bar{\gm}(\tau),\tau),\dot{\bar{\gm}}(\tau)).
\end{equation}Then $(\bar{\gm}(\tau),u(\tau),p(\tau))$  is $C^1$ and satisfies the characteristics equation (\ref{hjech}).

 By (L3) and Lemma \ref{tria1}, a standard argument (see \cite{F3,M1}) shows that the differentiability of $\bar{\gm}(\tau)$ for $\tau\in [t_0,t_0+\epsilon]$ can be extended to the whole interval $(0,t)$. So far, we complete the proof of Lemma \ref{chc}.\End

\subsection{Characteristics and fundamental solutions}
In this step, we will prove a relation between  $h_{x_0,u_0}(x,t)$ and $U(t)$, where  $U(t)$ belongs to a characteristic curve $(X(t),U(t),P(t))$. More precisely, we have the following lemma:
\begin{Lemma}\label{chh}
For $t\in (0, T]$, let $(X(t),U(t),P(t))$ denote a characteristic curve and $\mathcal{S}_{x_0,u_0}^x$ denote the set of $(X(t),U(t),P(t))$ satisfying $X(0)=x_0$, $X(t)=x$ and $U(0)=u_0$, then we have
\begin{equation}\label{hu11}
h_{x_0,u_0}(x,t)=\inf\left\{U(t):(X(t),U(t),P(t))\in \mathcal{S}_{x_0,u_0}^x\right\}.
\end{equation}
\end{Lemma}
\Proof
For the simplicity of notations, we use $C_i$ to denote the constants only depending on $t$. First of all, we prove that the infimum on the right side of (\ref{hu11}) can be achieved.
More precisely, there exists $(\bar{X}(t),\bar{U}(t),\bar{P}(t))$ such that
\begin{equation}
\bar{U}(t)=\inf\left\{U(t):(X(t),U(t),P(t))\in \mathcal{S}_{x_0,u_0}^x\right\}.
\end{equation}
In terms of the characteristic equation (\ref{hjech}), it follows from the Legendre transformation that
\begin{equation}\label{U1}
\dot{U}(t)=L(X(t),U(t),\dot{X}(t)).
\end{equation}
By the assumptions (L2) and (L4), there exists a constant $C_1$ such that for any $s\in [0,t]$, \[L(X(s),U(s),\dot{X}(s))\geq C_1,\]
 hence $\dot{U}(s)\geq C_1$. Moreover, a simple calculation implies
\begin{equation}
U(s)\geq u_0-|C_1|t.
\end{equation}
Therefore, $\inf_{\mathcal{S}_{x_0,u_0}^x}U(t)$ exists, which is denoted by $\tilde{u}$. Then, one can find a sequence $(X_n(t),U_n(t),\dot{X}_n(t))$ such that (extracting a subsequence if necessary)
$U_n(t)\rightarrow \tilde{u}$  as $n\rightarrow\infty$,
 hence, for $n$ large enough, we have
\begin{equation}\label{U2}
U_n(t)\leq \tilde{u}+1.
\end{equation}
From $\dot{U}(s)\geq C_1$, it follows that for any $s\in [0,t]$, \[U(t)-U(s)\geq C_1(t-s),\]
 which together with (\ref{U2}) implies
\begin{equation*}
U_n(s)\leq \tilde{u}+1+|C_1|t.
\end{equation*}It follows that for any $s\in [0,t]$,
$|U_n(s)|\leq C_2$.
Hence, according to (L2), it follows that for $n$ large enough, there exists $s_n\in [\frac{t}{2},t]$ such that
$|\dot{X}_n(s_n)|\leq C_3$.
Based on the compactness of $M\times [0,t]$, we have  (extracting a subsequence if necessary) as $n\rightarrow\infty$ \begin{equation}
(X_n(s_n),U_n(s_n),\dot{X}_n(s_n))\rightarrow (\bar{x},\bar{u},\dot{\bar{x}}).
\end{equation}By virtue of continuous dependence of solutions of (\ref{hjech}) on initial conditions, it follows that there exists a characteristic curve $(\bar{X}(t),\bar{U}(t),\dot{\bar{X}}(t))$ via Legendre transformation such that
\begin{equation}
\begin{array}{lll}
\bar{X}(0)=x_0,& \bar{X}(\bar{s})=\bar{x},& \bar{X}(t)=x,\\
\bar{U}(0)=u_0,& \bar{U}(\bar{s})=\bar{u},& \bar{U}(t)=\tilde{u}.\\
\end{array}
\end{equation}Therefore, there exists $(\bar{X}(t),\bar{U}(t),\bar{P}(t))$ such that
\begin{equation}\label{inf}
\bar{U}(t)=\inf\left\{U(t):(X(t),U(t),P(t))\in \mathcal{S}_{x_0,u_0}^x\right\}.
\end{equation}

In the following, we will prove
\begin{equation}
h_{x_0,u_0}(x,t)=\bar{U}(t).
\end{equation}

By virtue of (\ref{U1}), we have
\begin{equation}
\bar{U}(t)=u_0+\int_0^tL(\bar{X}(\tau),\bar{U}(\tau),\dot{\bar{X}}(\tau))d\tau.
\end{equation} By Lemma \ref{exist} and Lemma \ref{chc}, it follows that there exists a characteristic curve (via Legendre transformation) $(\tilde{X}(t),\tilde{U}(t),\dot{\tilde{X}}(t))$ such that  $\tilde{U}(t)=h_{x_0,u_0}(x,t)$ and
\begin{equation}\label{hxu}
h_{x_0,u_0}(x,t)=u_0+\int_0^tL(\tilde{X}(\tau),h_{x_0,u_0}(\tilde{X}(\tau),\tau),\dot{\tilde{X}}(\tau))d\tau.
\end{equation}
By contradiction, we assume $\bar{U}(t)\neq h_{x_0,u_0}(x,t)$. From (\ref{inf}), it follows that  $\bar{U}(t)<h_{x_0,u_0}(x,t)$. Based on Lemma \ref{chc}, for given $\hat{x}$, $\hat{u}$ and $\hat{t}$, there exist $\epsilon$, $\delta>0$ small enough such that for any $x\in \mathring{B}(\hat{x},\epsilon)$ and $s\in (\hat{t}-\delta,\hat{t}+\delta)$, we have
\begin{equation}\label{xus}
h_{\hat{x},\hat{u}}(x,s)=S(x,s),
\end{equation}where $S(x,t)$ is a classical solution of (\ref{hje}) satisfying $S(\hat{x},\hat{t})=\hat{u}$. Since  $(\bar{X}(t),\bar{U}(t),$ $\dot{\bar{X}}(t))$ is a $C^1$ characteristic curve, then it is easy to see that there exist $N>0$ and a partition as follows:
\begin{equation}
\{(\bar{X}(s_i),s_i):i=0,\ldots,N\},
\end{equation} such that $s_{i+1}\in (s_i-\delta,s_i+\delta)$ and $\bar{X}(s_{i+1})\in \mathring{B}(\bar{X}(s_i),\epsilon)$. In particular, we let $s_0=0$ and $s_N=t$. By (\ref{xus}), it yields
\begin{equation}
h_{\bar{X}(s_i),\bar{U}(s_i)}(\bar{X}(s_{i+1}),s_{i+1})=\bar{U}(s_{i+1}).
\end{equation}From
 (\ref{triin}), it follows that
\begin{align*}
h_{x_0,u_0}(x,t)-u_0&\leq \sum_{i=0}^N h_{\bar{X}(s_i),\bar{U}(s_i)}(\bar{X}(s_{i+1}),s_{i+1})-\bar{U}(s_i)\\
&=\bar{U}(t)-u_0,
\end{align*}which implies $h_{x_0,u_0}(x,t)\leq \bar{U}(t)$. It contradicts the assumption $\bar{U}(t)<h_{x_0,u_0}(x,t)$.
 This finishes the proof of Lemma \ref{chh}.
\End

So far, we have finished the proof of Theorem \ref{two}.

\section{\sc Representation of the viscosity solution}
In this section, we are devoted to proving Theorem \ref{one}. First of all, we construct a variational solution  of (\ref{hje}) with initial condition.
\subsection{Fundamental solutions and  variational solutions}

Based on Theorem \ref{two}, it follows that under the assumptions (L1)-(L4), there exists a unique $h_{y,\phi(y)}(x,t)\in C(M\times (0,T],\R)$ such that
\begin{equation}\label{hphi}
h_{y,\phi(y)}(x,t)=\phi(y)+\inf_{\substack{\gm(t)=x \\  \gm(0)=y} }\int_0^tL(\gm(\tau),h_{y,\phi(y)}(\gm(\tau),\tau),\dot{\gm}(\tau))d\tau,
\end{equation}where the infimums are taken among the absolutely continuous curves $\gm:[0,t]\rightarrow M$.

\begin{Lemma}\label{uh}
Let
\begin{equation}
u(x,t):=\inf_{y\in M}h_{y,\phi(y)}(x,t),
\end{equation}
then \begin{equation}\label{fixu}
u(x,t)=\inf_{\gm(t)=x}\left\{\phi(\gm(0))+\int_0^tL(\gm(\tau),u(\gm(\tau),\tau),\dot{\gm}(\tau))d\tau\right\}.
\end{equation}
\end{Lemma}

\Proof
We denote the Lipschitz constant by $\lambda$. The idea of the  proof is similar to the one in Lemma \ref{tria1}. By Lemma \ref{unique},   $\inf_{y\in M}h_{y,\phi(y)}(x,t)$ and $u(x,t)$ determined by (\ref{hphi}) and (\ref{fixu}) are unique. We will prove $\inf_{y\in M}h_{y,\phi(y)}(x,t)=u(x,t)$ in the following.

On one hand, we will prove
\[u(x,t)\geq\inf_{y\in M}h_{y,\phi(y)}(x,t).\]
 Let $\gm_1:[0,t]\rightarrow M$ be a calibrated curve of $u$ with $\gm_1(t)=x$. Let $\bar{y}:=\gm_1(0)$. It suffices to show
  \begin{equation}\label{s11uh}
u(x,t)\geq h_{\bar{y},\phi(\bar{y})}(x,t).
\end{equation}
By contradiction, we assume $u(x,t)< h_{\bar{y},\phi(\bar{y})}(x,t)$. By (\ref{hphi}) and (\ref{fixu}), we have
\begin{align*}
u(x,t)=\phi(\bar{y})+\int_0^{t}L(\gm_1(\tau),u(\gm_1(\tau),\tau),\dot{\gm}_1(\tau))d\tau.
\end{align*}
\begin{align*}
h_{\bar{y},\phi(\bar{y})}(x,t)\leq\phi(\bar{y})+
\int_0^{t}L(\gm_1(\tau),h_{\bar{y},\phi(\bar{y})}(\gm_1(\tau),\tau),\dot{\gm}_1(\tau))d\tau.
\end{align*}
For any $\sigma\in [0,t]$, we denote $\bar{u}(\sigma):=u(\gm_1(\sigma),\sigma)$ and $\bar{h}(\sigma):=h_{\bar{y},\phi(\bar{y})}(\gm_1(\sigma),\sigma)$. In particular, we have $\bar{u}(t)=u(x,t)$ and $\bar{h}(t)=h_{\bar{y},\phi(\bar{y})}(x,t)$.  Let
 \begin{equation}
F(\sigma):=\bar{h}(\sigma)-\bar{u}(\sigma),
\end{equation}where $\sigma\in [0,t]$. It is easy to see that $\bar{u}(0)=\phi(\bar{y})=\bar{h}(0)$. Then we have  $F(0)=0$. The assumption $u(x,t)< h_{\bar{y},\phi(\bar{y})}(x,t)$ implies $F(t)>0$. Hence, there exists $\sigma_0\in [0,t)$ such that $F(\sigma_0)=0$ and  $F(\sigma)> 0$ for $\sigma> \sigma_0$.
Moreover, for any $\tau\in (\sigma_0,t]$, we have
  \begin{equation}\label{s77}
\bar{u}(\tau)=\bar{u}(\sigma_0)+\int_{\sigma_0}^{\tau}L(\gm_1(\sigma),\bar{u}(\sigma),\dot{\gm}_1(\sigma))d\sigma.
\end{equation}
Let $\gm_2$ be a calibrated curve of $h_{\bar{y},\phi(\bar{y})}$ with $\gm_2(0)=\bar{y}$, $\gm_2(\sigma_0)=\gm_1(\sigma_0)$. For $\sigma\in [\sigma_0,\tau]$, we construct $\gm_\tau:[0,\tau]\rightarrow M$ as follows:
\begin{equation}
\gm_\tau(\sigma)=\left\{\begin{array}{ll}
\hspace{-0.4em}\gm_2(\sigma),& \sigma\in [0,\sigma_0],\\
\hspace{-0.4em}\gm_1(\sigma),&\sigma\in (\sigma_0,\tau].\\
\end{array}\right.
\end{equation}
Moreover, for any $\tau\in (\sigma_0,t]$, we have
  \begin{equation}
\bar{h}(\tau)\leq \bar{h}(\sigma_0)+\int_{\sigma_0}^{\tau}L(\gm_1(\sigma),\bar{h}(\sigma),\dot{\gm}_1(\sigma))d\sigma.
\end{equation}
Since $\bar{h}(\sigma_0)-\bar{u}(\sigma_0)=F(\sigma_0)=0$, a direct calculation implies
  \begin{equation}
\bar{h}(\tau)-\bar{u}(\tau)\leq \int_{\sigma_0}^\tau\lambda(\bar{h}(\sigma)-\bar{u}(\sigma))d\sigma.
\end{equation}
Hence, we have
  \begin{equation}
F(\tau)\leq \int_{\sigma_0}^\tau\lambda F(\sigma)d\sigma,
\end{equation}
which is
 in contradiction with Lemma \ref{gron}. Hence, we have
 \begin{equation}\label{s22uh}
u(x,t)\geq\inf_{y\in M}h_{y,\phi(y)}(x,t).
\end{equation}

On the other hand, we will prove $u(x,t)\leq\inf_{y\in M}h_{y,\phi(y)}(x,t)$. Due to the compactness of $M$, there exists $\tilde{y}$ such that $\inf_{y\in M}h_{y,\phi(y)}(x,t)=h_{\tilde{y},\phi(\tilde{y})}(x,t)$.
 Let $\gm_2:[0,t]\rightarrow M$ be a calibrated curve of $h$ with $\gm_2(0)=\tilde{y}$ and $\gm_2(t)=x$. Moreover, we have
 \begin{align*}
u(x,t)\leq\phi(\tilde{y})+\int_0^{t}L(\gm_2(\tau),u(\gm_2(\tau),\tau),\dot{\gm}_2(\tau))d\tau.
\end{align*}
\begin{align*}
h_{\bar{y},\phi(\tilde{y})}(x,t)=\phi(\tilde{y})+
\int_0^{t}L(\gm_2(\tau),h_{\tilde{y},\phi(\tilde{y})}(\gm_2(\tau),\tau),\dot{\gm}_2(\tau))d\tau.
\end{align*}
By a similar argument as (\ref{s77})-(\ref{s22uh}), we have
 \begin{equation}
u(x,t)\leq\inf_{y\in M}h_{y,\phi(y)}(x,t),
\end{equation}
which together with (\ref{s22uh}) implies
\[u(x,t)=\inf_{y\in M}h_{y, \phi(y)}(x,t).\]This finishes the proof of Lemma \ref{uh}.
\End

\begin{Lemma}\label{uisw}
$u(x,t)$ determined by (\ref{fixu}) is a variational solution of (\ref{hje}) with initial condition.
\end{Lemma}

\Proof
Let $\gm:[t_1,t_2]\rightarrow M$ be a continuous and piecewise $C^1$ curve and Let $\bar{\gm}:[0,t_1]\rightarrow M$ be a calibrated curve of $u$ satisfying $\bar{\gm}(t_1)=\gm(t_1)$. We construct a curve $\xi:[0,t_2]\rightarrow M$ defined as follows:
\begin{equation}\label{xi}
\xi(t)=\left\{\begin{array}{ll}
\hspace{-0.4em}\bar{\gm}(t),& t\in [0,t_1],\\
\hspace{-0.4em}\gm(t),&t\in (t_1,t_2].\\
\end{array}\right.
\end{equation}
From (\ref{fixu}), it follows that
\begin{align*}
u(\gm(t_2)&,t_2)-u(\gm(t_1),t_1)\\
&= \inf_{\gm_2(t_2)=\gm(t_2)}\left\{\phi(\gm_2(0))+
\int_0^{t_2}L(\gm_2(\tau),u(\gm_2(\tau),\tau),\dot{\gm}_2(\tau))d\tau\right\}\\
&\ \ \ \ -\inf_{\gm_1(t_1)=\gm(t_1)}\left\{\phi(\gm_1(0))+
\int_0^{t_1}L(\gm_1(\tau),u(\gm_1(\tau),\tau),\dot{\gm}_1(\tau))d\tau\right\},\\
&\leq \phi(\xi(0))+
\int_0^{t_2}L(\xi(\tau),u(\xi(\tau),\tau),\dot{\xi}(\tau))d\tau\\
&\ \ \ \ -\phi(\bar{\gm}(0))-
\int_0^{t_1}L(\bar{\gm}(\tau),u(\bar{\gm}(\tau),\tau),\dot{\bar{\gm}}(\tau))d\tau,
\end{align*}which together with (\ref{xi}) gives rise to
\begin{equation}
u(\gm(t_2),t_2)-u(\gm(t_1),t_1)\leq \int_{t_1}^{t_2}L(\gm(\tau),u(\gm(\tau),\tau),\dot{\gm}(\tau))d\tau,
\end{equation}which verifies (i) of Definition \ref{nw}.

By means of Lemma \ref{chc}, there exists a $C^1$ calibrated curve $\gm:[t_1,t_2]\rightarrow M$ with $\gm(t_2)=x$ such that
  \begin{equation}
u(x,t_2)-u(\gm(t_1),t_1)=\int_{t_1}^{t_2}L(\gm(\tau),u(\gm(\tau),\tau),\dot{\gm}(\tau))d\tau.
\end{equation}which implies (ii) of Definition \ref{nw}. This completes the proof of Lemma\ref{uisw}.
\End

\subsection{Variational solutions  and viscosity solutions}
In this subsection, we will prove the following lemma:
\begin{Lemma}\label{awisv}
A variational solution of (\ref{hje}) with initial condition is a viscosity solution.
\end{Lemma}
\Proof
 Let $u$ be a variational solution. Since $u(x,0)=\phi(x)$, then it suffices to consider $t\in (0,T]$.  We use $V\subset M$ to denote an open subset. Let $\phi:V\times\R\rightarrow \R$ be a $C^1$ test function such that $u-\phi$ has a maximum at $(x_0,t_0)$. This means $\phi(x_0,t_0)-\phi(x,t)\leq u(x_0,t_0)-u(x,t)$. Fix $v\in T_{x_0}M$ and for a given $\delta>0$, we choose a $C^1$ curve $\gm:[t_0-\delta,t_0+\delta]\rightarrow M$ with $\gm(t_0)=x_0$ and $\dot{\gm}(t_0)=\xi$. For $t\in [t_0-\delta,t_0]$, we have
\begin{align*}
\phi(\gm(t_0),t_0)-\phi(\gm(t),t)&\leq u(\gm(t_0),t_0)-u(\gm(t),t),\\
&\leq \int_t^{t_0}L(\gm(\tau),u(\gm(\tau),\tau),\dot{\gm}(\tau))d\tau,
\end{align*}where the second inequality is based (i) of Definition \ref{nw}. Hence,
\begin{equation}
\frac{\phi(\gm(t),t)-\phi(\gm(t_0),t_0)}{t-t_0}\leq\frac{1}{t-t_0}
\int_{t_0}^tL(\gm(\tau),u(\gm(\tau),\tau),\dot{\gm}(\tau))d\tau.
\end{equation}Let $t\rightarrow t_0$, we have
\[\partial_t\phi(x_0,t_0)+\partial_{x}\phi(x_0,t_0)\cdot \xi\leq L(x_0,u(x_0,t_0),\xi),\]which together with Legendre transformation implies
\[\partial_t\phi(x_0,t_0)+H(x_0,u(x_0,t_0),\partial_x\phi(x_0,t_0))\leq 0,\]which shows that $u$ is a viscosity subsolution.

To complete the proof of Theorem \ref{awisv}, it remains to show that $u$ is a supersolution. $\psi:V\times\R\rightarrow \R$ be a $C^1$ test function and $u-\psi$ has a minimum at $(x_0,t_0)$. We have $\psi(x_0,t_0)-\psi(x,t)\geq u(x_0,t_0)-u(x,t)$. From (ii) of Definition \ref{nw}, there exists a $C^1$ curve $\gm:[0,t_0]\rightarrow M$ with $\gm(t_0)=x_0$ and $\dot{\gm}(t_0)=\eta$ such that for $0\leq t<t_0$, we have
  \begin{equation}
u(\gm(t_0),t_0)-u(\gm(t),t)=\int_{t}^{t_0}L(\gm(\tau),u(\gm(\tau),\tau),\dot{\gm}(\tau))d\tau.
\end{equation}Hence
\[\psi(x_0,t_0)-\psi(x,t)\geq\int_{t}^{t_0}L(\gm(\tau),u(\gm(\tau),\tau),\dot{\gm}(\tau))d\tau.\]
Moreover, we have
\[\frac{\psi(\gm(t),t)-\psi(\gm(t_0),t_0)}{t-t_0}\geq\frac{1}{t-t_0}
\int_{t_0}^tL(\gm(\tau),u(\gm(\tau),\tau),\dot{\gm}(\tau))d\tau.\]Let $t$ tend to $t_0$, it gives rise to
\[\partial_t\psi(x_0,t_0)+\partial_{x}\psi(x_0,t_0)\cdot \eta\geq L(x_0,u(x_0,t_0),\eta),\]which implies
\[\partial_t\phi(x_0,t_0)+H(x_0,u(x_0,t_0),\partial_x\phi(x_0,t_0))\geq 0.\]This finishes the proof of Lemma \ref{awisv}.
\End

By the comparison theorem (see \cite{Ba3} for instance), it yields that  the viscosity solution of (\ref{hje}) is unique under the assumptions (H1)-(H4). So far, we have obtained that there exists a unique  viscosity solution $u(x,t)$ of (\ref{hje}) with initial condition $u(x,0)=\phi(x)$.
So far,  we complete the proof of Theorem \ref{one}.

\section{\sc Solution semigroup}
In this section, we will prove Theorem \ref{four}. Let $u(x,t)$ be the unique viscosity solution of (\ref{hje}) with initial condition $u(x,0)=\phi(x)$. We introduce an implicitly defined nonlinear operator $T_t$ such that
\begin{equation}\label{uttt}
u(x,t)=T_t\phi(x).
\end{equation}
It follows from (\ref{fixu}) that
\begin{equation}\label{51}
T_t\phi(x)=\inf_{\gm(t)=x}\left\{\phi(\gm(0))+\int_0^tL(\gm(\tau),T_\tau\phi(\gm(\tau)),
\dot{\gm}(\tau))d\tau\right\}.
\end{equation}
where the infimums are taken among absolutely continuous curves. In particular, the infimums are attained at the characteristics of (\ref{hje}).
The following lemma implies $T_t$ is a semigroup.
\begin{Proposition}\label{semi}
$\{T_t\}_{t\geq 0}$ is a one-parameter semigroup of operators from
$C(M,\R)$ into itself.
\end{Proposition}
\Proof
It is easy to see $T_0  = Id$. It suffices to prove that
$T_{t+s}=T_t\circ T_s$ for any $t,s\geq 0$.

For every $\eta\in C(M, \mathbb{R})$ and $u \in C(M\times
[0,T],\mathbb{R})$, we define an operator $\A_t^{\eta}$ such that
\begin{equation}
\A^{\eta}[u](x,t) = \inf_{\gamma(t)=x} \left\{ \eta\big(\gamma(0)\big) + \int_0^t
L\big(\gamma(\tau),u\big(\gamma(\tau),\tau\big), \dot\gamma(\tau)\big) d\tau
\right\}.
\end{equation}
By virtue of Theorem \ref{one}, it follows that $\A^{\eta}$ has a unique fixed point.

By (\ref{51}), we have
\begin{equation*}
\begin{split}
T_t\circ T_s \phi(x) &= \inf_{\gamma(t)=x} \left\{
T_s\phi\big({\gamma}(0)\big) + \int_0^{t}
L\big({\gamma}(\tau),T_{\tau}\circ T_s
\phi\big({\gamma}(\tau)\big), \dot{{\gamma}}(\tau)\big)d\tau
\right\}\\
& = \A^{T_s \phi} [T_t\circ T_s \phi](x).
\end{split}
\end{equation*}

On the other hand,
\begin{equation*}
\begin{split}
T_{t+s} \phi(x) &= \inf_{\substack{\gamma(t+s)=x}} \left\{ \phi\big(\gamma(0)\big) +
\int_0^{t+s} L\big(\gamma(\tau),T_\tau
\phi\big(\gamma(\tau)\big),
\dot\gamma(\tau)\big)d\tau \right\}\\
& = \inf_{\substack{\gamma(t+s)=x}}
\left\{ \phi\big(\gamma(0)\big) + \bigg(\int_0^{s} +
 \int_{s}^{t+s}\bigg) L\big(\gamma(\tau),T_\tau
\phi\big(\gamma(\tau)\big), \dot\gamma(\tau)\big) d\tau \right\}\\
& = \inf_{\substack{\gamma(t+s)=x}}
\left\{ T_s\phi\big(\gamma(s)\big) + \int_s^{t+s}
L\big(\gamma(\tau),T_\tau \phi\big(\gamma(\tau)\big),
\dot\gamma(\tau)\big)d\tau \right\}\\
& = \inf_{\substack{\bar{\gamma}(t)=x}} \left\{ T_s\phi\big(\bar{\gamma}(0)\big) +
\int_0^{t} L\big(\bar{\gamma}(\tau),T_{\tau+s}
\phi\big(\bar{\gamma}(\tau)\big),
\dot{\bar{\gamma}}(\tau)\big)d\tau \right\}  \\
& = \A^{T_s \phi} [T_{t+s} \phi](x).
\end{split}
\end{equation*}
Hence, both $T_t\circ T_s \phi$ and $T_{t+s} \phi$ are fixed points of $\A^{T_s \phi} $, which together with the uniqueness of the fixed point of $\A^{T_s \phi} $ yields  $T_{t+s}=T_t\circ T_s$. This completes
the proof of  Proposition \ref{semi}.
\End

\begin{Proposition}\label{htth}
For given $y\in M$ and $s>0$, we have
\[h_{y,\phi(y)}(x,s+t)=T_th_{y,\phi(y)}(x,s),\]
where $h_{y,\phi(y)}(\cdot,\cdot)$ is defined as (\ref{hphi}).
\end{Proposition}
\Proof
Fix $y\in M$ and $s>0$, one can define a continuous function $h_{y,\phi(y)}^s: M\rightarrow \R$ by $h_{y,\phi(y)}^s(x)=h_{y,\phi(y)}(x,s)$. Based on Lemma \ref{uh} and (\ref{uttt}), we have
\[T_th_{y,\phi(y)}^s(x)=\inf_{z\in M}h_{z,h_{y,\phi(y)}^s(z)}(x,t).\]
It follows from Lemma \ref{tria1} that
\[h_{y,\phi(y)}(x,s+t)=\inf_{z\in M}h_{z,h_{y,\phi(y)}^s(z)}(x,t).\]
Hence, we have
\[h_{y,\phi(y)}(x,s+t)=T_th_{y,\phi(y)}(x,s).\]
This completes the proof of Proposition \ref{htth}.
\End

To fix the notion, we call $T_t$ a solution semigroup. In the following subsections, we will prove some further properties of the solution semigroup $T_t$.

First of all, it is easy to obtain the following proposition about the monotonicity of $T_t$.
\begin{Proposition}[Monotonicity]\label{monotonicity}
For given $\phi,\psi\in C(M,\R)$ and $t\geq 0$, if $\phi\leq\psi$, then $T_t\phi\leq T_t\psi$.
\end{Proposition}
\Proof
 For given $\phi,~\psi\in C(M,\mathbb{R})$ with
$\phi\leq \psi$, by contradiction, we assume that there exist
$t_1>0$ and $x_1\in M$ such that $T_{t_1}\phi(x_1) > T_{t_1}
\psi(x_1)$.
Let $\gm_\psi: [0,t_1]\rightarrow M$ be a calibrated curve of $T_t\psi$ with $\gm_\psi(t_1)=x_1$. We denote
\[F(\tau)=T_{\tau}\phi(\gm_\psi(\tau)) - T_{\tau}
\psi(\gm_\psi(\tau)).\] It is easy to see that $F(\tau)$ is continuous and $F(t_1)>0$. Since
 \[F(0)=\phi(\gm_\psi(0)) -
\psi(\gm_\psi(0))\leq 0,\]
there exists $t_0\in [0,t_1)$ such that $F(t_0)=0$ and for any $\tau\in [t_0,t_1]$, $F(\tau)\geq 0$.
It follows from (\ref{51}) that for any $s\in [t_0,t_1]$, we have
\begin{equation}\label{monotonicity estimate}
\begin{split}
&T_{s}\phi(\gm_\psi(s)) - T_{s}\psi(\gm_\psi(s)) \\= &\inf_{\substack{\gamma(s)=\gm_\psi(s)}} \left\{ T_{t_0}\phi(\gamma({t_0}))
+ \int_{t_0}^{s} L(\gamma(\tau),T_\tau
\phi(\gamma(\tau)),
\dot\gamma(\tau)) d\tau \right\} \\
&-\inf_{\substack{\gamma(s)=\gm_\psi(s)\\}}
\left\{ T_{t_0}\psi(\gamma({t_0})) + \int_{t_0}^{s}
L(\gamma(\tau),T_\tau
\psi(\gamma(\tau)), \dot\gamma(\tau))  d\tau \right\} \\
 \leq &  T_{t_0}\phi(\gm_\psi({t_0})) -  T_{t_0}\psi(\gm_\psi({t_0})) \ +\\
&\int_{t_0}^{s} L(\gm_\psi(\tau),T_\tau
\phi(\gm_\psi(\tau)), \dot\gm_\psi(\tau)) -
L(\gm_\psi(\tau),T_\tau \psi(\gm_\psi(\tau)),
\dot{\gm}_\psi(\tau)) d\tau,\\
\leq &\int_{t_0}^{s} L(\gm_\psi(\tau),T_\tau
\phi(\gm_\psi(\tau)), \dot\gm_\psi(\tau)) -
L(\gm_\psi(\tau),T_\tau \psi(\gm_\psi(\tau)),
\dot{\gm}_\psi(\tau)) d\tau,\\
\leq & \lambda \int_{t_0}^{s}T_\tau
\phi(\gm_\psi(\tau))-T_\tau \psi(\gm_\psi(\tau))d\tau,
\end{split}
\end{equation}
where $\lambda$ denotes  the Lipschitz constant of $L$.

Hence, we have
\begin{equation}
F(s)\leq\lambda\int_{t_0}^{s}F(\tau)d\tau,
\end{equation}
which is in contradiction with Lemma \ref{gron}. Therefore, we have
\[T_{t_1}\phi(x_1) \leq T_{t_1}\psi(x_1).\]
This finishes the proof of Proposition \ref{monotonicity}.\End

Proposition \ref{monotonicity} can be viewed as a comparison principle for (\ref{hje}).
By a similar argument as the one in Proposition \ref{monotonicity}, one can obtain the Lipschitz continuity of $T_t$. For $\phi\in C(M,\R)$, we use $\|\phi\|_\infty$ to denote $C^0$-norm of $\phi$. We have the following proposition.
\begin{Proposition}[Lipschitz continuity]\label{nonex}
For given $\phi,\psi\in C(M,\R)$ and $t\geq 0$, we have $\|T_t\phi-T_t\psi\|_\infty\leq e^{\lambda t}\|\phi-\psi\|_\infty$.
\end{Proposition}

First of all, we prove the following lemma.
\begin{Lemma}\label{hinequ}
For given $(x,t)\in M\times (0,T]$, $x_0\in M$, $u,v\in\R$, we have
\begin{equation}
|h_{x_0,u}(x,t)-h_{x_0,v}(x,t)|\leq e^{\lambda t}|u-v|.
\end{equation}
\end{Lemma}
\Proof
A similar argument implies the monotonicity of $h_{x_0,u}(x,t)$ with respect to $u$. More precisely, if $u\geq v$, then $h_{x_0,u}(x,t)\geq h_{x_0,v}(x,t)$.

Since $u,v\in\R$, then we have the dichotomy: a) $u\leq v$, b) $u>v$. For Case a), we have $h_{x_0,u}(x,t)\leq h_{x_0,v}(x,t)$. Let $\gm_u$ be a calibrated curve of $h_{x_0,u}$ with $\gm_u(0)=x_0$ and $\gm_u(t)=x$. From the monotonicity of $h_{x_0,u}(x,t)$, it follows that for any $s\in (0,t]$,
\begin{equation}\label{tpro}
h_{x_0,u}(\gm_u(s),s)\leq h_{x_0,v}(\gm_u(s),s).
\end{equation}
In terms of the definition of $h_{x_0,u}(x,t)$, we have
\begin{align*}
&h_{x_0,v}(\gm_u(s),s)- h_{x_0,u}(\gm_u(s),s)\\
\leq &v-u+\int_0^sL(\gm_u(\tau),h_{x_0,v}(\gm_u(\tau),\tau),\dot{\gm}_u(\tau))
-L(\gm_u(\tau),h_{x_0,u}(\gm_u(\tau),\tau),\dot{\gm}_u(\tau))d\tau,\\
\leq & v-u+\int_0^s\lambda|h_{x_0,v}(\gm_u(\tau),\tau)- h_{x_0,u}(\gm_u(\tau),\tau)|d\tau.
\end{align*}
Let $F(\tau):=h_{x_0,v}(\gm_u(\tau),\tau)- h_{x_0,u}(\gm_u(\tau),\tau)$. It follows from (\ref{tpro}) that $F(\tau)\geq 0$ for any $\tau\in (0,t]$. Hence, we have
\[F(s)\leq v-u+\int_0^s\lambda F(\tau)d\tau.\]
By Gronwall's inequality, it yields
\begin{equation}
F(s)\leq (v-u)e^{\lambda s}.
\end{equation}
In particular, we verify Lemma \ref{hinequ} for Case a).

For Case b), we have $h_{x_0,u}(x,t)\geq h_{x_0,v}(x,t)$. Let $\gm_v$ be a calibrated curve of $h_{x_0,v}$ with $\gm_v(0)=x_0$ and $\gm_v(t)=x$. Let $G(\tau):=h_{x_0,u}(\gm_u(\tau),\tau)- h_{x_0,v}(\gm_u(\tau),\tau)$. By a similar argument as Case a), we have
\begin{equation}
G(s)\leq (u-v)e^{\lambda s}.
\end{equation}

Therefore, we completes the proof of Lemma \ref{hinequ} for any $u,v\in\R$.
\End

\noindent \textbf{Proof of Proposition \ref{nonex}: }
By Theorem \ref{one}, we have
\[u(x,t)=\inf_{y\in M}h_{y,\phi(y)}(x,t),\]
where $u(x,t)$ is a viscosity solution of (\ref{hje}).
It follows from (\ref{uttt}) that $u(x,t)=T_t\phi(x)$. Hence, we have
\begin{equation}\label{thphi}
T_t\phi(x)=\inf_{y\in M}h_{y,\phi(y)}(x,t).
\end{equation}
Similarly, we have
\begin{equation}\label{thpsi}
T_t\psi(x)=\inf_{z\in M}h_{z,\psi(z)}(x,t).
\end{equation}
Lemma \ref{hinequ} implies that $h_{y,\phi(y)}(x,t)$ is continuous with respect to $y$. Based on the compactness of $M$, the infimums in (\ref{thphi}) and (\ref{thpsi}) can be attained at $y_0$ and $z_0$ respectively. On one hand, we have
\begin{align*}
&T_t\phi(x)-T_t\psi(x)\\
\leq &h_{z_0,\phi(z_0)}(x,t)-h_{z_0,\psi(z_0)}(x,t),\\
\leq & e^{\lambda t}|\phi(z_0)-\psi(z_0)|,\\
\leq & e^{\lambda t}\|\phi(x)-\psi(x)\|_{\infty}.
\end{align*}
On the other hand, we have
\begin{align*}
&T_t\phi(x)-T_t\psi(x)\\
\geq &h_{y_0,\phi(y_0)}(x,t)-h_{y_0,\psi(y_0)}(x,t),\\
\geq & -e^{\lambda t}|\phi(y_0)-\psi(y_0)|,\\
\geq & -e^{\lambda t}\|\phi(x)-\psi(x)\|_{\infty}.
\end{align*}
Hence,
\[\|T_t\phi(x)-T_t\psi(x)\|_{\infty}\leq e^{\lambda t}\|\phi(x)-\psi(x)\|_{\infty}.\]
This completes the proof of Proposition \ref{nonex}.
\End

So far, we have finished the proof of Theorem \ref{four}.
\section{\sc Large time behavior of the viscosity solution }
In this section, we will prove Theorem \ref{five}, which is concerned with the large time behavior of $T_t$.
\subsection{Critical values}
For $c\in\R$, we denote $L_c:=L+c$. For given $x_0, u_0, x, t$ where $t\in (0,+\infty)$, we denote
\begin{equation}
h_{x_0,u_0}^c(x,t)=u_0+\inf_{\substack{\gm(t)=x \\  \gm(0)=x_0} }\int_0^tL_c(\gm(\tau),h_{x_0,u_0}^c(\gm(\tau),\tau),\dot{\gm}(\tau))d\tau,
\end{equation}where the infimums are taken among the absolutely continuous curves $\gm:[0,t]\rightarrow M$.
The critical value set is defined as
\begin{equation}
\mathcal{C}=\left\{c: |h_{x_0,u_0}^c(x,t)|\leq K(u_0)\right\},
\end{equation}
where  $K(u_0)$ is a positive constant depending on $u_0$.
For $a\in\R$, we use $c(L(x,a,\dot{x}))$ to denote Ma\~{n}\'{e} critical value of $L(x,a,\dot{x})$. By \cite{SWY}, $\mathcal{C}\neq\emptyset$ if $L(x,u,\dot{x})$ is non-increasing with respect to $u$. The following theorem gives the more general conditions under which $\mathcal{C}$ is not empty.
\begin{Theorem}\label{ubound}
If there exist $L_1(x,\dot{x})$, $L_1(x,\dot{x})$ satisfying (L1) and (L2) such that
\begin{equation}
\begin{cases}
\lim_{a\rightarrow+\infty}L(x,a,\dot{x})\leq L_1(x,\dot{x}),\\
\lim_{a\rightarrow-\infty}L(x,a,\dot{x})\geq L_2(x,\dot{x}),\\
c(L_1(x,\dot{x}))\geq c(L_2(x,\dot{x})),
\end{cases}
\end{equation}
then for any $c\in [c(L_2(x,\dot{x})),c(L_1(x,\dot{x}))]$, there exists $K>0$ such that for $t\geq\delta$
\begin{equation}
 |h_{x_0,u_0}^c(x,t)-u_0|\leq K.
\end{equation}
\end{Theorem}
\Proof
For the simplicity of notations, we denote $c_1:=c(L_1(x,\dot{x}))$ and $c_2:=c(L_2(x,\dot{x}))$.  The proof is divided into two steps.

In Step One, we prove $h_{x_0,u_0}^c(x,t)$ is upper bounded.
 By contradiction, we assume that for any $K>0$, there exists $t'\gg 0$   such that $h_{x_0,u_0}^c(x,t')=K+u_0$. Let $\gm:[0,t']\rightarrow M$ be a calibrated curve of $h_{x_0,u_0}^c(x,t) $ with $\gm(0)=x_0$ and $\gm(t')=x$. Since $h_{x_0,u_0}^c(\gm(t),t)=U(t)$ is $C^1$ with respect to $t$ for $t\in (0,t')$, then one  can find $t''>\delta$ such that $h_{x_0,u_0}^c(\gm(t''),t'')=K/2+u_0$ and for any $\tau\in [t'',t']$,
\begin{equation}\label{ulb}
\frac{K}{2}\leq h_{x_0,u_0}^c(\gm(\tau),\tau)-u_0\leq K.
\end{equation}

Based on the definition of $h^{x_0,u_0}_c(x,t)$, we have
\begin{equation}
h_{x_0,u_0}^c(x,t')=u_0+\inf_{\substack{\gm(t')=x \\  \gm(0)=x_0} }\int_0^{t'}L_c(\gm(\tau),h_{x_0,u_0}^c(\gm(\tau),\tau),\dot{\gm}(\tau))d\tau.
\end{equation}
Let $\gm_1:[t'',t']\rightarrow M$ be a minimal curve with $\gm_1(t'')=\gm(t'')$ and $\gm_1(t')=x$. That is
\begin{equation}
\int_0^{t'}L_1(\gm_1(\tau),\dot{\gm}_1(\tau))d\tau=\inf_{\substack{\tilde{\gm}(t')=x \\  \tilde{\gm}(t'')=\gm(t'')} }\int_0^{t'}L_1(\tilde{\gm}(\tau),\dot{\tilde{\gm}}(\tau))d\tau.
\end{equation}

Moreover, it follows that
\begin{align*}
h_{x_0,u_0}^c(x,t')&=h_{x_0,u_0}^c(\gm(t''),t'')+\int_{t''}^{t'}
L_c(\gm(\tau),h_{x_0,u_0}^c(\gm(\tau),\tau),\dot{\gm}(\tau))d\tau,\\
&\leq \frac{K}{2}+u_0+\int_{t''}^{t'}
L_c(\gm_1(\tau),h_{x_0,u_0}^c(\gm_1(\tau),\tau),\dot{\gm}_1(\tau))d\tau,\\
&\leq \frac{K}{2}+u_0+\int_{t''}^{t'}
L_1(\gm_1(\tau),\dot{\gm}_1(\tau))+cd\tau,\\
&\leq \frac{K}{2}+u_0+\int_{t''}^{t'}
L_1(\gm_1(\tau),\dot{\gm}_1(\tau))+c_1d\tau,\\
&=\frac{K}{2}+u_0+h_{c_1}^{t'-t''}(\gm(t''),x),
\end{align*}
where the third inequality is owing to (\ref{ulb}) and the assumption
\[\lim_{a\rightarrow+\infty}L(x,a,\dot{x})\leq L_1(x,\dot{x}).\]
$h_{c_1}^{t'-t''}(\gm(t''),x)$ denotes the minimal action with respect to $L_1$ (see \cite{F3} for instance). It is easy to see that $t'-t''>1$, then from the compactness of $M$, it follows that $h_{c_1}^{t'-t''}(\gm(t''),x)$ has a bound denoted by $A$ independent of $t',t''$ and $x$. Hence, we have
\[K\leq \frac{K}{2}+A.\]
Since $K$ is large enough, then we have a contradiction if we take $K>3A$.

In Step Two, we prove $h_{x_0,u_0}^c(x,t)$ is lower bounded.
By contradiction, we assume that for any $-K<0$, there exists $t'\gg 0$   such that $h_{x_0,u_0}^c(x,t')=-K+u_0$. Let $\gm:[0,t']\rightarrow M$ be a calibrated curve of $h_{x_0,u_0}^c(x,t) $ with $\gm(0)=x_0$ and $\gm(t')=x$. Hence, one  can find $t''>\delta$ such that $h_{x_0,u_0}^c(\gm(t''),t'')=-K/2+u_0$ and for any $\tau\in [t'',t']$,
\begin{equation}\label{ulb1}
-K\leq h_{x_0,u_0}^c(\gm(\tau),\tau)-u_0\leq -\frac{K}{2}.
\end{equation}

Moreover, it follows that
\begin{align*}
h_{x_0,u_0}^c(x,t')&=h_{x_0,u_0}^c(\gm(t''),t'')+\int_{t''}^{t'}
L_c(\gm(\tau),h_{x_0,u_0}^c(\gm(\tau),\tau),\dot{\gm}(\tau))d\tau,\\
&\geq -\frac{K}{2}+u_0+\int_{t''}^{t'}
L_2(\gm(\tau),\dot{\gm}(\tau))+cd\tau,\\
&\geq -\frac{K}{2}+u_0+\int_{t''}^{t'}
L_2(\gm(\tau),\dot{\gm}(\tau))+c_2d\tau,\\
&=-\frac{K}{2}+u_0+h_{c_2}^{t'-t''}(\gm(t''),x),
\end{align*}
where the second inequality is owing to (\ref{ulb1}) and the assumption
\[\lim_{a\rightarrow-\infty}L(x,a,\dot{x})\geq L_2(x,\dot{x}).\]
$h_{c_2}^{t'-t''}(\gm(t''),x)$ denotes the minimal action with respect to $L_2$. It is easy to see that $h_{c_2}^{t'-t''}(\gm(t''),x)>0$ (see \cite{F3}). Hence, we have
\[-K\geq -\frac{K}{2},\]
which contradicts the assumption $K>0$. So far, we have shown $h_{x_0,u_0}^c(x,t)$ is uniformly bounded for $t\geq \delta$.
This finishes the proof of Theorem \ref{ubound}.
\End

Unfortunately, we do not know whether the conditions in Theorem \ref{ubound} is sharp or not. Generally, the critical value may not exist. For instance, we consider the Hamilton-Jacobi equation:
\begin{equation}
\begin{cases}
\partial_tu-u+\frac{1}{2}|\partial_xu|^2=c,\\
u(x,0)=\phi(x),
\end{cases}
\end{equation}
where the Hamiltonian $H(x,u,p)=\frac{1}{2}|p|^2-u$ which satisfies the assumptions (H1)-(H4). It is easy to see that for any non-constant function $\phi(x)\in C(M,\R)$ and $c\in \R$, there holds $u(x,t)\rightarrow\infty$ exponentially as $t\rightarrow\infty$, which means $\mathcal{C}=\emptyset$. Conversely, we are concerned with
\begin{equation}
\begin{cases}
\partial_tu+u+\frac{1}{2}|\partial_xu|^2=c,\\
u(x,0)=\phi(x).
\end{cases}
\end{equation}
It follows that for any $\phi(x)\in C(M,\R)$ and $c\in \R$, there holds $u(x,t)\rightarrow c$ as $t\rightarrow\infty$, for which $\mathcal{C}=\R$.
 Moreover, $u(x,t)$ converges to $u(x)\equiv c$ which is the unique viscosity solution of the stationary equation on $M$:
\[u+\frac{1}{2}|\partial_xu|^2=c.\]

\subsection{Large time behavior of the solution semigroup}
In order to consider the large time behavior of $T_t$, we need the following assumption:
\begin{itemize}
\item [\textbf{(H5)}]   \textbf{Non-emptiness}: The critical value set $\mathcal{C}$ is not empty.
\end{itemize}
 Without ambiguity, we still use $L$ instead of $L_c$ to denote $L+c$ for $c\in \mathcal{C}$. The same to $H$ and $T_t$. Based on a similar argument as the one in \cite{SWY}, we have
\begin{Proposition}\label{unilip}
For any $\phi(x)\in C(M,\R)$, $x,y\in M$ and $t> \delta$, we have
\begin{itemize}
\item [I.] there exists a positive constant $K$ independent of $t$ such that $\|h_{y,\phi(y)}(x,t)\|_\infty\leq K(\phi)$;
\item [II.] for $\delta>0$, there exists a compact subset $\mathcal{K}_\delta$ such that for every calibrated curve $\gm$ of $h_{y,\phi(y)}(x,t)$ and any $t>\delta$, we have
\begin{equation*}
(\gm(t),h_{y,\phi(y)}(\gm(t),t),\dot{\gm}(t))\in \mathcal{K}_\delta.
\end{equation*}
\item [III.] for  $\delta>0$, the family of functions $(x,t)\rightarrow h_{y,\phi(y)}(x,t)$ is equi-Lipschitz on $(x,t)\in M\times [\delta,+\infty)$.
\end{itemize}
\end{Proposition}
We omit the proof of Proposition \ref{unilip} here. See Theorem 1.3 in \cite{SWY} for the details. By Theorem \ref{four}, there holds
\[T_t\phi(x)=\inf_{y\in M}h_{y,\phi(y)}(x,t).\]
Let
\[\underline{u}(x):=\liminf_{t\rightarrow\infty}T_t\phi(x).\]
Proposition \ref{unilip} implies $\underline{u}(x)$ is a Lipschitz function.

\begin{Lemma}\label{fixpoint}
For any $t\geq 0$, we have
\[T_t\underline{u}(x)=\underline{u}(x).\]
\end{Lemma}
\Proof
We denote $u_s(x):=\inf_{t\geq s}T_t\phi(x)$, then $\underline{u}(x)=\lim_{s\rightarrow\infty}u_s(x)$. It suffices to prove $T_\delta\underline{u}(x)=\underline{u}(x)$ for any $\delta\geq 0$.

One one hand, we will prove $T_\delta\underline{u}(x)\leq\underline{u}(x)$. It is easy to see that
\begin{equation}\label{61}
\begin{split}
u_{\delta+s}(x)&=\inf_{t\geq \delta+s}T_t\phi(x),\\
&=\inf_{t-s\geq \delta}T_\delta\circ T_{t-\delta}\phi(x),\\
&=\inf_{t\geq s}T_\delta\circ T_t\phi(x).
\end{split}
\end{equation}
For $t\geq s$, we have $u_s(x)\leq T_t\phi(x)$. It follows from the monotonicity of $T_t$ that
\[T_\delta\circ T_t\phi(x)\geq T_\delta u_s(x).\]
Moreover, we have
\[\inf_{t\geq s}T_\delta\circ T_t\phi(x)\geq T_\delta u_s(x),\]
which together with (\ref{61}) implies
\[u_{\delta+s}(x)\geq T_\delta u_s(x).\]
 Taking the limit as $s\rightarrow\infty$ in both sides, we have
\begin{equation}\label{65}
T_\delta\underline{u}(x)\leq\underline{u}(x).
\end{equation}

On the other hand, we have
\begin{equation}\label{62}
T_\delta u_s(x)=T_\delta \left(\inf_{t\geq s}T_t\phi(x)\right)=\inf_{y\in M}h_{y,\inf_{t\geq s}T_t\phi(y)}(x,\delta).
\end{equation}

\textbf{Claim:}
\[\inf_{y\in M}h_{y,\inf_{t\geq s}T_t\phi(y)}(x,\delta)=\inf_{t\geq s}\left(\inf_{y\in M}h_{y,T_t\phi(y)}(x,\delta)\right).\]

\Proof
It is easy to see that $h_{y,\inf_{t\geq s}T_t\phi(y)}(x,\delta)$ is continuous with respect to $y$. Based on the compactness of $M$, there exists $y_0\in M$ such that the infimum is attained. Hence,
\[\inf_{y\in M}h_{y,\inf_{t\geq s}T_t\phi(y)}(x,\delta)=h_{y_0,\inf_{t\geq s}T_t\phi(y_0)}(x,\delta)\geq\inf_{t\geq s}\left(\inf_{y\in M}h_{y,T_t\phi(y)}(x,\delta)\right).\]

On the other hand, it follows from the monotonicity of $h_{y,u}(x,t)$ with respect to $u$ that
\[h_{y,\inf_{t\geq s}T_t\phi(y)}(x,\delta)\leq h_{y,T_t\phi(y)}(x,\delta),\]
which yields
\[\inf_{y\in M}h_{y,\inf_{t\geq s}T_t\phi(y)}(x,\delta)\leq\inf_{t\geq s}\left(\inf_{y\in M}h_{y,T_t\phi(y)}(x,\delta)\right).\]
This completes the proof of the claim.
\End

From (\ref{62}), we have
\begin{align*}
T_\delta u_s(x)&=\inf_{y\in M}h_{y,\inf_{t\geq s}T_t\phi(y)}(x,\delta)=\inf_{t\geq s}\left(\inf_{y\in M}h_{y,T_t\phi(y)}(x,\delta)\right),\\
&= \inf_{t\geq s}T_\delta\circ T_t\phi(x)=\inf_{t\geq \delta+s}T_t\phi(x),\\
&=u_{\delta+s}(x)\geq u_s(x).
\end{align*}
Taking the limit as $s\rightarrow\infty$ in both sides, it follows that
\begin{equation}\label{64}
T_\delta\underline{u}(x)\geq\underline{u}(x),
\end{equation}
 which together with (\ref{65}) completes the proof of Lemma \ref{fixpoint}.
\End

\begin{Lemma}\label{fixpointweak}
$T_tu(x)=u(x)$ for any $t\geq 0$ if and only if $u(x)$ is a weak KAM solution of the following stationary equation:
\begin{equation}\label{stahh}
H(x, u(x), \partial_x u(x))=0.
\end{equation}
\end{Lemma}
\Proof
We suppose $T_tu(x)=u(x)$ for any $t\geq 0$. By virtue of a similar argument as Lemma \ref{uisw},  it yields that for each continuous piecewise $C^1$ curve $\gm:[t_1,t_2]\rightarrow M$ where $0\leq t_1<t_2\leq T$, we have
\begin{equation}
u(\gm(t_2))-u(\gm(t_1))\leq\int_{t_1}^{t_2}L(\gm(\tau),
u(\gm(\tau)),\dot{\gm}(\tau))d\tau,
\end{equation}which implies (i) of Definition \ref{weakkam}. In addition,
there exists a $C^1$ calibrated curve $\gm_t:[-t,0]\rightarrow M$ with $\gm_t(0)=x$ such that for any $t'\in [-t,0]$, we have
  \begin{equation}
u(x)-u(\gm_t(t'))=\int_{t'}^{0}L(\gm_t(\tau),u(\gm_t(\tau)),\dot{\gm}_t(\tau))d\tau.
\end{equation}
Based on the a priori compactness given by Lemma \ref{unilip} II., for a given $\delta>0$, there exists a compact subset $\mathcal{K}_\delta$ such that for  any $s>\delta$, we have
\begin{equation*}
(\gm_t(s),u(\gm_t(s)),\dot{\gm}_t(s))\in \mathcal{K}_\delta.
\end{equation*}
Since $\gm_t$ is a calibrated curve, it follows from Lemma \ref{chc} that
\[(\gm_t(s),u(\gm_t(s)),\dot{\gm}_t(s))=\Phi_s(\gm_t(0),u(\gm_t(0)),\dot{\gm}_t(0))=\Phi_s(x,u(x),\dot{\gm}_t(0)).\]
The points $(\gm_t(0),u(\gm_t(0)),\dot{\gm}_t(0))$ are contained in a compact subset, then one can find a sequence $t_n$ such that $(x,\dot{\gm}_{t_n}(0))$ tends to $(x,v_\infty)$ as $n\rightarrow \infty$. Fixing $t'\in (-\infty,0]$, the function $s\mapsto \Phi_s(x,u(x),\dot{\gm}_{t_n}(0))$ is defined on $[t',0]$ for $n$ large enough. By the continuity of $\Phi_s$, the sequence converges uniformly on the compact interval $[t',0]$ to the map $s\mapsto \Phi_s(x,v_\infty)$. Let
\[(\gm_\infty(s),u(\gm_\infty(s)), \dot{\gm}_\infty(s)):=\Phi_s(x,v_\infty),\] then
for any $t'\in (-\infty,0]$, we have
  \begin{equation}
u(x)-u(\gm_\infty(t'))=\int_{t'}^{0}L(\gm_\infty(\tau),u(\gm_\infty(\tau)),\dot{\gm}_\infty(\tau))d\tau,
\end{equation}
which implies (ii) of Definition \ref{weakkam}.
Hence, $u$ is a weak KAM solution of (\ref{stahh}).

Conversely, we suppose $u$ is a weak KAM solution of (\ref{stahh}). By (i) of Definition \ref{weakkam}, we have
\[u\leq T_tu.\]
 By (ii) of Definition \ref{weakkam}, for any $x\in M$, there exists a $C^1$ curve $\bar{\gm}:(-\infty,0]\rightarrow M$ with $\bar{\gm}(0)=x$ such that for any $t\in [0,+\infty)$,
  \begin{equation}
u(x)-u(\bar{\gm}(-t))=\int^{0}_{-t}L(\bar{\gm}(\tau),u(\bar{\gm}(\tau)),\dot{\bar{\gm}}(\tau))d\tau.
\end{equation}
We define the curve $\gm:[0,t]\rightarrow M$ by $\gm(s)=\bar{\gm}(s-t)$. There hold $\gm(t)=x$ and
\[u(x)=\int^{t}_{0}L(\gm(\tau),u(\gm(\tau)),\dot{\gm}(\tau))d\tau+u(\gm(0)),\]
 which implies
 \[T_tu\leq u.\]
This completes the proof of Lemma \ref{fixpointweak}.
\End

So far, we have finished the proof of Theorem \ref{five}.
\section{\sc Projected Aubry set and the stationary equation}
In this section, we will define the projected Aubry set with respect to the stationary Hamilton-Jacobi equation (\ref{station}). Moreover, we will prove Theorem \ref{six}. Without ambiguity, we still use $L$ instead of $L_c$ to denote $L+c$ for $c\in \mathcal{C}$. The same to $H$ and $T_t$.
\subsection{Projected Aubry set}
For a given $s_0>0$, we take $\phi(x)=h_{x_0,u_0}(x,s_0)$, where $x_0\in M, u_0\in \R$. Then $\phi(x)\in C(M,\R)$. By virtue of Theorem \ref{four}, for $t\geq 0$, there holds
\[h_{x_0,u_0}(x,s_0+t)=T_th_{x_0,u_0}(x,s_0)=T_t\phi(x).\]
Let \[h_{x_0,u_0}(x,\infty):=\liminf_{t\rightarrow\infty}h_{x_0,u_0}(x,t).\] It is easy to see that
\[h_{x_0,u_0}(x,\infty)=\liminf_{t\rightarrow\infty}h_{x_0,u_0}(x,s_0+t)=\liminf_{t\rightarrow\infty}T_t\phi(x).\]
By Theorem \ref{five}, $h_{x_0,u_0}(x,\infty)$ is a viscosity solution of (\ref{station}). Based on \cite{SWY}, ``liminf" can be replaced with ``lim" if $H(x,u,p)$ is non-decreasing with respect to $u$.

We denote
 \[B(x,u;y):=h_{x,u}(y,\infty)-u.\]

 $B(x,u;y)$ can be referred as the barrier function dented by $h^\infty(x,y)$ in Mather-Fathi theory. Moreover, we define the projected Aubry set as follows;
\begin{equation*}
\mathcal{A}:=\{(x,u)\in M\times\R\ |\ B(x,u;x)=0\}.
\end{equation*}
\begin{Proposition}\label{anoe}
$\mathcal{A}\neq\emptyset$  under the assumptions (H1)-(H5).
\end{Proposition}
\Proof
For given $x_0,y_0\in M$ and $u_0\in \R$, let $\gm_n:[0,t_n]\rightarrow M$ be a calibrated curve of $h_{x_0,u_0}$ with $\gm_n(0)=x_0$ and $\gm_n(t_n)=y_0$.
Under the assumptions (H1)-(H5), Proposition \ref{unilip} implies that for $\delta>0$, there exists a compact subset $\mathcal{K}_\delta\subset M\times\R$ such that for any $\tau>\delta$, we have
\begin{equation*}
(\gm_n(\tau),h_{x_0,u_0}(\gm_n(\tau),\tau))\in \mathcal{K}_\delta.
\end{equation*}
Extracting a subsequence if necessary, one can find $s_n,\tau_n\in [0,t_n]$ satisfying $s_n-\tau_n\rightarrow\infty$ as $n\rightarrow\infty$ and
\[d\left((\gm_n(s_n),h_{x_0,u_0}(\gm_n(s_n),s_n)),(\gm_n(\tau_n),h_{x_0,u_0}(\gm_n(\tau_n),\tau_n))\right)\leq \frac{1}{n},\]
where $d(\cdot,\cdot)$ denotes the distance induced by a Riemannian metric on $M\times\R$.
Hence, there exists $(x_\infty,u_\infty)\in M\times\R$ such that
\[h_{x_\infty,u_\infty}(x_\infty,\infty)=u_\infty.\]
This completes the proof of Proposition \ref{anoe}.
\End

\subsection{Projected Aubry set and  the viscosity solution}
In order to find finer properties of $\mathcal{A}$, we add the following assumption:
\begin{itemize}
\item [\textbf{(H6)}]  \textbf{Strict increase}: $H(x,u,p)$ is strictly increasing with respect to $u$ for a given $(x,p)\in T^*M$.
\end{itemize}
It is easy to see that (H5) holds under (H6) (see \cite{SWY}). (H6) is equivalent to
\begin{itemize}
\item [\textbf{(L6)}]  \textbf{Strict decrease}: $L(x,u,\dot{x})$ is strictly decreasing with respect to $u$ for a given $(x,\dot{x})\in TM$.
\end{itemize}

%
%
%

It is easy to obtain the following proposition about the  contractibility of $T_t$. For $\phi\in C(M,\R)$, we use $\|\phi\|_\infty$ to denote $C^0$-norm of $\phi$. We have the following proposition.
\begin{Lemma}[Contractibility]\label{nonexst}
For given $\phi,\psi\in C(M,\R)$ and $t\geq 0$, if there exists $\bar{x}$ such that $\phi(\bar{x})\neq \psi(\bar{x})$, then  $\|T_t\phi-T_t\psi\|_\infty<\|\phi-\psi\|_\infty$.
\end{Lemma}

First of all, we prove the following lemma.
\begin{Lemma}\label{hinequst}
For given  $x_0,x\in M$, $u,v\in\R$ and $t>0$, if $u\neq v$, then we have
\begin{equation}
|h_{x_0,u}(x,t)-h_{x_0,v}(x,t)|<|u-v|.
\end{equation}
\end{Lemma}
\Proof
A similar argument as Proposition \ref{monotonicity} implies the monotonicity of $h_{x_0,u}(x,t)$ with respect to $u$. More precisely, if $u\geq v$, then $h_{x_0,u}(x,t)\geq h_{x_0,v}(x,t)$.

Since $u,v\in\R$ and $u\neq v$, then we have the dichotomy: a) $u< v$, b) $u>v$. For Case a), we have $h_{x_0,u}(x,t)\leq h_{x_0,v}(x,t)$. Let $\gm_u$ be a calibrated curve of $h_{x_0,u}$ with $\gm_u(0)=x_0$ and $\gm_u(t)=x$. From the monotonicity of $h_{x_0,u}(x,t)$ with respect to $u$ and the continuity of $h_{x_0,u}(x,t)$ with respect to $t$, it follows that there exits $\delta>0$ such that for any $s\in (0,\delta]$,
\begin{equation}\label{tprost}
h_{x_0,u}(\gm_u(s),s)< h_{x_0,v}(\gm_u(s),s).
\end{equation}
In terms of the definition of $h_{x_0,u}(x,t)$, we have
 \begin{align*}
 h_{x_0,v}(x,t)&=v+\inf_{\substack{\gm(t)=x \\  \gm(0)=x_0} }\int_0^tL(\gm(\tau),h_{x_0,v}(\gm(\tau),\tau),\dot{\gm}(\tau))d\tau,\\
 &\leq v+\int_0^tL(\gm_u(\tau),h_{x_0,v}(\gm_u(\tau),\tau),\dot{\gm}_u(\tau))d\tau,\\
 &<v+\int_0^tL(\gm_u(\tau),h_{x_0,u}(\gm_u(\tau),\tau),\dot{\gm}_u(\tau))d\tau,\\
 &=v-u+u+\int_0^tL(\gm_u(\tau),h_{x_0,u}(\gm_u(\tau),\tau),\dot{\gm}_u(\tau))d\tau,\\
 &=v-u+ h_{x_0,u}(x,t),
 \end{align*}
 where the third inequality is from (L6) and (\ref{tprost}). Hence, there holds
 \begin{equation}
0\leq  h_{x_0,v}(x,t)- h_{x_0,u}(x,t)<v-u.
\end{equation}
In particular, we verify Lemma \ref{hinequst} for Case a).

For Case b), we have $h_{x_0,u}(x,t)\geq h_{x_0,v}(x,t)$. Let $\gm_v$ be a calibrated curve of $h_{x_0,v}$ with $\gm_v(0)=x_0$ and $\gm_v(t)=x$. It follows that there exits $\delta'>0$ such that for any $s\in (0,\delta']$,
\begin{equation}\label{tprost}
h_{x_0,v}(\gm_v(s),s)< h_{x_0,u}(\gm_v(s),s).
\end{equation} By a similar argument as Case a), we have
 \begin{equation}
0\leq  h_{x_0,u}(x,t)- h_{x_0,v}(x,t)<u-v.
\end{equation}

Therefore, we completes the proof of Lemma \ref{hinequst} for any $u,v\in\R$ and $u\neq v$.
\End

\noindent \textbf{Proof of Lemma \ref{nonexst}: }
By Theorem \ref{one} and Theorem \ref{four}, we have
\begin{equation}\label{thphist}
T_t\phi(x)=\inf_{y\in M}h_{y,\phi(y)}(x,t),\quad
T_t\psi(x)=\inf_{z\in M}h_{z,\psi(z)}(x,t).
\end{equation}
 Based on the compactness of $M$, the infimums in (\ref{thphist})  can be attained at $y_0$ and $z_0$ respectively. If $\phi(z_0)=\psi(z_0)$, we have
\begin{align*}
&T_t\phi(x)-T_t\psi(x)\\
\leq &h_{z_0,\phi(z_0)}(x,t)-h_{z_0,\psi(z_0)}(x,t)=  0,\\
< &|\phi(\bar{x})-\psi(\bar{x})|
\leq \|\phi(x)-\psi(x)\|_{\infty}.
\end{align*}
 If $\phi(z_0)\neq\psi(z_0)$, we have
 \begin{align*}
&T_t\phi(x)-T_t\psi(x)\\
\leq &h_{z_0,\phi(z_0)}(x,t)-h_{z_0,\psi(z_0)}(x,t),\\
< &|\phi(z_0)-\psi(z_0)|
\leq \|\phi(x)-\psi(x)\|_{\infty}.
\end{align*}
Hence, there holds
\begin{equation}\label{tstt}
T_t\phi(x)-T_t\psi(x)<\|\phi(x)-\psi(x)\|_{\infty}.
\end{equation}
Similarly, we have
\begin{equation}
T_t\phi(x)-T_t\psi(x)>-\|\phi(x)-\psi(x)\|_{\infty},
\end{equation}
which together with (\ref{tstt}) implies
\[\|T_t\phi(x)-T_t\psi(x)\|_{\infty}<\|\phi(x)-\psi(x)\|_{\infty}.\]
This completes the proof of Lemma \ref{nonexst}.
\End
\begin{Lemma}[Uniqueness]\label{stuni}
Let $H$ satisfy (H1)-(H4) and (H6), then there exists a unique viscosity solution satisfying the stationary equation
\begin{equation}\label{hxup}
H(x,u(x),\partial_xu(x))=0.
\end{equation}
\end{Lemma}
\Proof
A similar argument as Theorem 7.6.2 in \cite{F3} implies the weak KAM solution is also equivalent to the viscosity solution for the stationary equation (\ref{hxup}). By Lemma \ref{fixpointweak}, $u$ is a weak KAM solution if and only if $T_tu=u$ for any $t\geq 0$.

By contradiction, we assume there exist at least two weak KAM solutions $u(x),v(x)\in C(M,\R)$ with $u(x_0)\neq v(x_0)$ for some $x_0\in M$. Hence, there hold $T_tu=u$ and $T_tv=v$. It follows from Proposition \ref{nonexst} that
\[\|u-v\|_{\infty}=\|T_tu-T_tv\|_{\infty}<\|u-v\|_{\infty},\]
which is a contradiction. This completes the proof of Lemma \ref{stuni}.
\End

We use $\pi:M\times\R\rightarrow M$ to denote the standard projection via $(x,u)\rightarrow x$.
\begin{Lemma}\label{aubry}
Let $H$ satisfy (H1)-(H4) and (H6),  for $x\in \pi\mathcal{A}$, there exists a unique $u_x$ such that $(x,u_x)\in \mathcal{A}$. Moreover, there holds $u(x)=u_x$ for any $x\in \pi\mathcal{A}$, where $u(x)\in C(M,\R)$ is the unique viscosity solution.
\end{Lemma}
\Proof
By contradiction, we assume for $x_0\in \pi\mathcal{A}$, there exist $u_1,u_2\in \R$ such that $h_{x_0,u_1}(x_0,\infty)=u_1$ and $h_{x_0,u_2}(x_0,\infty)=u_2$. By Theorem \ref{five}, $h_{x_0,u_1}(x,\infty)$ and
$h_{x_0,u_2}(x,\infty)$ are viscosity solutions of (\ref{hxup}). Lemma \ref{stuni} implies
\[h_{x_0,u_1}(x,\infty)\equiv h_{x_0,u_2}(x,\infty).\]
In particular, we have
\[u_1=h_{x_0,u_1}(x_0,\infty)= h_{x_0,u_2}(x_0,\infty)=u_2.\]

For $(y,u_{y}),(z,u_{z})\in \mathcal{A}$, we denote
\[u_y(x):=h_{y,u_{y}}(x,\infty),\quad u_z(x):= h_{z,u_z}(x,\infty).\]
It follows from Lemma \ref{stuni} that $u_y(x)\equiv u_z(x)$ denoted by $u(x)$. Then $u(x)=u_x$ for any $x\in \pi\mathcal{A}$.
\End

Based on Lemma \ref{aubry}, one can obtain the following lemma.
\begin{Lemma}[Representation formula]\label{repre}
Let $u(x)$ be the viscosity solution of (\ref{hxup}), then
\begin{equation}
u(x)=\inf_{y\in \pi\mathcal{A}}h_{y,u(y)}(x,\infty).
\end{equation}
\end{Lemma}
\Proof
By virtue of the uniqueness of the viscosity solution, it follows that for given $(x_0,u_0)\in M\times\R$, $u(x)=h_{x_0,u_0}(x,\infty)$.  Lemma \ref{tria1} yields
  \begin{equation}
h(x,t+s)=\inf_{y\in M}h_{y,h(y,t)}(x,s),
\end{equation}
where the ``inf" can be attained at $y$ belonging to the calibrated curve $\gm$ of $h$ with $\gm(0)=x_0$ and $\gm(t+s)=x$. Hence, we have
\begin{equation}\label{111}
u(x)=h_{x_0,u_0}(x,\infty)=\inf_{y\in M}h_{y,h_{x_0,u_0}(y,\infty)}(x,\infty).
\end{equation}

On one hand, it follows from Lemma \ref{aubry} that
\begin{equation}\label{222}
\inf_{y\in M}h_{y,h_{x_0,u_0}(y,\infty)}(x,\infty)\leq \inf_{y\in \pi\mathcal{A}}h_{y,h_{x_0,u_0}(y,\infty)}(x,\infty)=\inf_{y\in \pi\mathcal{A}}h_{y,u(y))}(x,\infty).
\end{equation}
On the other hand, there exists $y_0\in \pi\mathcal{A}$ such that
\begin{align*}
\inf_{y\in M}h_{y,h_{x_0,u_0}(y,\infty)}(x,\infty)&=h_{y_0,h_{x_0,u_0}(y_0,\infty)}(x,\infty)=h_{y_0,u(y_0)}(x,\infty),\\
&\geq \inf_{y\in \pi\mathcal{A}}h_{y,u(y))}(x,\infty),
\end{align*}
which together with (\ref{111}) and (\ref{222}) yields
\[u(x)=\inf_{y\in \pi\mathcal{A}}h_{y,u(y)}(x,\infty).\]This completes the proof of Lemma \ref{repre}.
\End

So far, we have finished the proof of Theorem \ref{six}.

 \vspace{2ex}
\noindent\textbf{Acknowledgement}
This work is partially under the support of National Natural Science Foundation of China (Grant No. 11171071,  11325103) and
National Basic Research Program of China (Grant No. 11171146).

\vspace{2em}

{\sc Lin Wang}

{\sc School of Mathematical Sciences, Fudan University,
Shanghai 200433,
China.}

 {\it E-mail address:} \texttt{linwang.math@gmail.com}

\vspace{1em}

{\sc Jun Yan}

{\sc School of Mathematical Sciences, Fudan University,
Shanghai 200433,
China.}

 {\it E-mail address:} \texttt{yanjun@fudan.edu.cn}

\end{document}